\documentclass[twoside,11pt]{article}

\usepackage{blindtext}

 \usepackage[abbrvbib, preprint]{jmlr2e}

\usepackage{lastpage}
\jmlrheading{??}{??}{1-\pageref{LastPage}}{??; Revised ??}{??}{21-0000}{Marie Du Roy de Chaumaray, Salima El Kolei, Marie-Pierre Etienne and Matthieu Marbac}

\ShortHeadings{Estimation of the Order of Non-Parametric Hidden Markov Models }{Du Roy, El Kolei, Etienne and Marbac}
\firstpageno{1}

\usepackage[utf8]{inputenc}
\usepackage{natbib}
\usepackage[english]{babel}
\usepackage{fullpage}

\usepackage{amsfonts,amssymb,amsmath,amscd,latexsym,dsfont}
\usepackage{graphicx}
\usepackage{natbib}
\usepackage{color}
\usepackage{enumitem}
\usepackage{url}

\newcommand{\bA}{\boldsymbol{A}}
\newcommand{\bV}{\boldsymbol{V}}
\newcommand{\bZ}{\boldsymbol{Z}}
\newcommand{\bv}{\boldsymbol{v}}
\newcommand{\bw}{\boldsymbol{w}}

\newcommand{\bW}{\boldsymbol{W}}
\newcommand{\bM}{\boldsymbol{M}}
\newcommand{\bY}{\boldsymbol{Y}}
\newcommand{\bX}{\boldsymbol{X}}
\newcommand{\bx}{\boldsymbol{x}}
\newcommand{\by}{\boldsymbol{y}}
\newcommand{\ba}{\boldsymbol{a}}
\newcommand{\bb}{\boldsymbol{b}}
\newcommand{\bz}{\boldsymbol{z}}
\newcommand{\bu}{\boldsymbol{u}}
\newcommand{\bpi}{\boldsymbol{\pi}}
\newcommand{\bc}{\boldsymbol{c}}
\newcommand{\wwr}{\widehat{r}}

\newcommand{\spa}{L^2(\mathbb{R}^d)}
\newcommand{\pkg}{\texttt{HMMselect}}
\newcommand{\dd}{\, \mathrm{d}}
\newtheorem{assumption}{Assumption}
 \DeclareMathOperator*{\argmin}{arg\,min}

\begin{document}

\title{Estimation of the Order of Non-Parametric Hidden Markov Models using the Singular Values of an Integral Operator}

\author{\name Marie Du Roy de Chaumaray \email marie.du-roy-de-chaumaray@univ-rennes2.fr\\
       \addr Mathematical Research Institute of Rennes IRMAR\\
        Rennes University, Rennes, France
       \AND
       \name Salima El Kolei \email salima.el-kolei@ensai.fr\\
       \addr Univ.  Rennes, Ensai, CNRS, CREST\\
        UMR 9194, F-35000 Rennes, France        
 \AND
       \name Marie-Pierre Etienne \email marie-pierre.etienne@institut-agro.fr \\
       \addr Mathematical Research Institute of Rennes IRMAR\\
       Rennes University, Rennes, France
       \AND
       \name Matthieu Marbac \email matthieu.marbac-lourdelle@ensai.fr\\
       \addr Univ.  Rennes, Ensai, CNRS, CREST\\
       UMR 9194, F-35000 Rennes, France
       }

\editor{My editor}

\maketitle

\begin{abstract}
 Interested in estimating the order of a finite-state Hidden Markov Model (HMM) with nonparametric emission distributions, a new method that only requires full rank transition matrix and not linearly dependent emission distributions is introduced. This   method  relies on the equality between the order of the HMM and the rank of a specific integral operator. Since only the empirical counter-part of the singular values of the operator can be obtained, a thresholding procedure is proposed. At a non-asymptotic level, an upper-bound on the probability of overestimating the order of the HMM is provided. At an asymptotic level,  the consistency of the estimator is established. In addition a general heuristic that can be successfully applied to several problems in spectral analysis for designing a data-driven procedure for the threshold is introduced.  The approach has the advantage of not requiring any  knowledge of an upper-bound on the order of the HMM. Moreover, different types of data (including circular or mixed-type data) can be managed. The relevance of the approach is illustrated on numerical experiments and on real data considering multivariate data with directional variables.
\end{abstract}


\section{Introduction}
A discrete-time homogeneous hidden Markov model (HMM) defines the distribution of an observed process $(\bY_t)_{t\in\mathbb{N}}$ and a latent process $(X_t)_{t\in\mathbb{N}}$ such that the sequence of unobserved states $(X_t)_{t\in\mathbb{N}}$ follows a Markov chain and the observations $(\bY_t)_{t\in\mathbb{N}}$ are independent given the state sequence $(X_t)_{t\in\mathbb{N}}.$ The conditional distribution of $Y_t$, called emission distribution, only depend on the current state $X_t$. 
This  paper  focuses  on  finite  state  HMMs  where  the  latent process has a finite state space $\{1,\ldots,L\}$, the integer $L$ being called the order of the HMM. In this framework, the model is completely described by the order $L$, the initial distribution and the transition matrix of the hidden chain, and the emission distributions. Since the marginal distribution of each $\bY_t$ is a finite mixture model, finite state HMMs can be seen as an extension of finite mixture models where the assumption of independence between observations is relaxed (\emph{i.e.,} $\bY_t$ and $\bY_{t'}$ are not independent).
HMMs are a popular tool for modeling the dependency structure for univariate and multivariate processes driven by a latent Markov chain (see \citet{juang1991hidden,yang1995application,krogh2001predicting,choo2004recent,zucchini2009hidden} for examples of applications). They also provide tractable models for circular time series widely used in biology, meteorology and climate applications to model for instance the speed and the direction of wind, ocean current, or animal movements (see \citet{holzmann2006hidden, bulla2012multivariate, mastrantonio2016hidden}). 
Inferring the right order of the latent chain is an important issue, which precedes the estimation of the model parameters and their interpretation.
This paper focuses on the estimation of the order $L$ from univariate and multivariate data $(\bY_t)_{t\in\mathbb{N}}$ in a non-parametric setting. The estimation of $L$ is indeed achieved without any parametric assumption on the emission distributions, since we only require the linear independence between their probability distribution functions.

Initial developments on HMMs have been made in a parametric framework, which considers that the emission distributions belong to some given parametric distribution family. Considering the order of the HMM as known, the existing literature provides the parameter identifiability \citep{petrie1969probabilistic}, the algorithm for assessing the maximum likelihood estimator (MLE; \citet{baum1970maximization}), the consistency of the MLE \citep{leroux1992consistent} and its asymptotic normality  \citep{bickel1998asymptotic}. The identification of the order is more challenging and represents a difficult task, mainly because of a loss of identifiability of the model parameters when the order is overestimated.  The standard assumptions used to control the likelihood ratio test statistics,  are thus not satisfied when the order is overestimated. For instance,  \citet{gassiat2000likelihood} show that this statistic can diverge even for bounded parameters. Note that this issue already appears when estimating the number of components in parametric finite mixture models \citep{ciuperca2002likelihood}. Therefore, the order of parametric HMMs can be estimated by homogeneous tests \citep{holzmann2016testing}, penalized likelihood approaches \citep{volant2014hidden} or cross-validation approaches \citep{celeux2008selecting}. Using tools from information theory, \citet{gassiat2003optimal} have shown the strong consistency of the estimator of the order obtained by penalized maximum likelihood. Moreover, Bayesian approaches can be used by penalizing the likelihood and thus avoiding the issues due to the lack of identifiability of the parameters when the order is overestimated \citep{gassiat2014posterior}. Alternatively, \citet{robert2000bayesian} propose a Bayesian inference of the order $L$ through a reversible jump Markov Chain Monte Carlo method (MCMC). In \citet{chopin2007inference}, the author also proposes  a Bayesian strategy based on sequential Monte Carlo filter and MCMC. 
However, all these approaches consider parametric emission distributions but it is not always possible to restrict the model to such a convenient finite-dimensional space. Moreover they provide biased results when their parametric assumptions are violated. In such cases, non-parametric approaches can be used to model the emission distributions.

Non-parametric HMMs have been proved to be useful in a wide range of applications (see \citet{zhao2011nonparametric} for financial applications, \citet{couvreur2000wavelet} for voice activity detection, \citet{lambert2003non} for climate state identification and \citet{yau2011bayesian} for genomic applications). Nevertheless, identifiability of the parameters of finite state HMMs with non-parametric emission distributions has been investigated recently. \citet{gassiat2016nonparametric} consider the case of translation HMMs. They show that all the model parameters (including the infinite dimensional parameters) are identifiable as soon as the matrix that defines the joint distribution of two consecutive latent variables, is non-singular and the translation parameters are distinct. Note that their conditions are weaker than those used to obtain identifiability for location-scale mixture models. Indeed, for the latter, constraints must be added such as considering symmetric distributions \citep{hunter2007inference}. This additional assumption is no longer required for translation HMMs because of the dependency between a pair of consecutive observations. Based on the results of parameter identifiability for a mixture of products of univariate distributions \citep{AllmanAOS2009}, \citet{gassiat2016inference} state weaker sufficient conditions for parameter identifiability since they consider a full rank transition matrix of the latent chain and linearly independent emission probability distributions. 
The method introduced by the present paper for estimating the order of an HMM, is developed under these assumptions. Note that the assumptions made on the emission distributions have been weakened again by \citet{alexandrovich2016nonparametric} since they only require that the emission distributions are different. 

To estimate the (finite and infinite dimensional) parameters of non-parametric HMMs, kernel-based \citep{bonhomme2016estimating} or wavelet-based \citep{jin2006non} approaches can be used. 
Alternatively, \citet{bonhomme2016estimating} and \citet{de2017consistent} extended the spectral method proposed by \citet{hsu2012spectral} for estimating parametric HMMs, in order to deal with a non-parametric framework. However, all these methods are developed for a known order of the HMM. 
Estimating the order of a generic non-parametric HMM is still a challenging problem and to the best of our knowledge \citet{lehericy2019consistent} is the only paper to consider this problem in this non-parametric setting. The author proposes two methods that provide strongly consistent estimators of the order of the HMM. The first method considers a minimization of a penalized least-square criterion that relies on a projection of the emission distributions onto a family of nested parametric subspaces. For each subspace and each number of latent states, the criterion used for model selection is computed by minimizing the empirical counterpart of the penalized $L^2$ distance. Thus, the method provides an estimator of the order of the HMM together with estimators of the emission distributions.  The second method uses an estimator of the rank of a matrix computed from the distribution of a pair of consecutive observations. More precisely, this method relies on a spectral approach applied on the matrix containing the coordinates of the density of a pair of consecutive observations in some orthonormal basis. Thus, this method could be seen as an extension of the spectral method described in the Section 5 of Supplementary material of \citet{BonhommeJRSSB2016} to HMM. 
These two methods are complementary in practice. Indeed, numerical experiments presented in \citet{lehericy2019consistent} show that the penalized least-square method is more efficient for moderate sample sizes. Indeed, the non-convex criterion raises many problems for the minimization, in practice. To overcome this difficulty, the author proposes to use an approximate minimization algorithm (see \citet{hansen2011cma}) that requires a good initial condition since it might otherwise remain trapped in a local minima which renders this approach time-consuming for multivariate data and large sample.
Furthermore, considering all the subspaces and all the possible numbers of latent states 
makes this method computationally greedy. Therefore, the spectral method should be considered for large sample sizes. Both methods involve an unknown tuning parameter (\emph{i.e.,} constant in the penalty term of the penalized criterion and threshold for the spectral method) but also choices of the subspaces (\emph{i.e.,} family of nested parametric subspaces or the orthonormal basis) that can highly impact the results (see our numerical experiments).

This paper introduces a new simple method for selecting the order of a non-parametric HMM, by using the rank of an integral operator relying on the distribution of a pair of consecutive observations. This approach is inspired by the one proposed in \citet{kwon2021estimation} to estimate the number of components in nonparametric i.i.d. mixture models but we go further in adapting this framework for dependent and latent observations leading to new and different theoretical results. Furthermore, we propose a more general heuristics for designing a data-driven procedure which can be successfully applied to several problems in spectral analysis and give theoretical guarantees.
The interest of this approach from integral operators lies in the fact that unlike most of the spectral methods based on noisy matrices \citep{bonhomme2016estimating,de2017consistent,lehericy2019consistent}, the method does not require any choices of a functional basis or its number of elements. Hence, the proposed method does not require any knowledge of an upper bound of the order of the HMM. Moreover, different types of data (including circular or mixed-type data) can be managed.  
Since the distribution of the pair of consecutive observations is estimated with kernel method, only the empirical counter-part of the singular values of the operator can be obtained, we propose to use our new data-driven method for the thresholding procedure.
At a non-asymptotic level, an upper-bound on the probability of overestimating the order of the HMM is provided. At an asymptotic level, the consistency of the estimator is established. The control at non-asymptotic and asymptotic levels are obtained by a concentration inequality of the Hilbert-Schmidt norm of the empirical version of the operator. The statistical tools needed to establish these results differ from those used in \citet{kwon2021estimation}, and consequently the results are different.
Thus, using concentration results specific to Markov chains, a concentration inequality is obtained by considering a sum of two terms, where one term does not depend on the bandwidth and the second term does not depend on the probability of overestimating the order. Note that the bound obtained in the i.i.d. context considers a product between the bandwidth, the probability of overestimating the order and the sample size. This bound contains only terms that depends on the kernel and the bandwidth. Contrary to this setting and because of the dependency of observations the concentration inequality that we obtain depends on some unknown constant of the HMM (\emph{e.g.,} the mixing time). 
To circumvent this issue and practical convenience, we propose a data-driven procedure based on an unsupervised classification of the singular values of the operator and computed on mini-batches, for estimating the constant in the concentration inequality. Note that in \cite{lehericy2019consistent} the model selection for the spectral method is also based on a thresholding rule applied on the singular values whose choice is a delicate issue since it depends on the functional basis and on the number of elements. Hence, in his paper the author proposes an empirical method based on slope heuristic for the practical application. However, this approach requires an additional tuning parameter that states the number of singular values used to apply the slope heuristic. In theory, for the spectral methods to work, the rank of the spectral matrix needs to be equal to the order of the chain. Thus, it is necessary that the number of elements of the orthonormal basis tends to infinity, otherwise we only obtain an estimator of an upper-bound of the order. However, defining the thresholding rule for the case of increasing number of basis elements is still an open problem for the spectral methods. Indeed, for instance, the rank study performed in \citet{Kleibergen2006generalized} should be extended to matrices with increasing dimension (but fixed rank). Thus, in practice, the number of basis elements is set a priori. This number corresponds to an upper-bound on the order of the HMM. To the best of our knowledge, since the proposed method avoids the use of functional basis, it is the first method which does not make assumptions on an upper-bound of the order to be estimated. 
Numerical studies illustrate the relevance of this proposal and show also
that this new data-driven procedure guarantees good results for our estimator, but also improves the spectral results of \cite{lehericy2019consistent}.

This paper is organized as follows. 
Section~\ref{sec:rank} introduces the specific integral operator. 
Section~\ref{sec:th} presents the finite-sample size and the asymptotic properties of the estimator (including its consistency). 
Section~\ref{sec:cst} describes the new data-driven procedure with a theoretical justification.
Section~\ref{sec:computation} is devoted to the computational aspects of the methods. 
Section~\ref{sec:ex} illustrates the consistency of the estimator on simulated data and shows the relevance of the proposed method on benchmark data (including circular data). 
Section~\ref{sec:appli} shows the contribution of our approach on one real-life data set. 
Section~\ref{sec:cl} gives a conclusion and all the proofs are given in Appendix.

\section{Order of a HMM and rank of integral operators} \label{sec:rank}
\subsection{Hidden Markov model}
Let $\bY=(\bY_1^\top,\ldots,\bY_{n+1}^\top)^\top$ be a stationary sequence of random vectors $\bY_t$, where $\bY_t\in\mathbb{R}^d$ follows a finite state hidden Markov model (HMM) with $L$ latent states. This model assumes that there exists a stationary Markov chain $\bX=(X_1,\ldots,X_{n+1})^\top$ that is unobserved, where $X_t\in\{1,\ldots,L\}$. Moreover, conditionally on $\bX$, the $\bY_t$'s are independent and their distribution only depends on the current state $X_t$. The Markov chain is defined by a full rank transition matrix $\bA$ having $\bpi=(\pi_1,\ldots,\pi_L)^\top$ as stationary distribution. Finally, the densities of the emission distributions $f_1,\ldots,f_L$ are assumed to be linearly independent, where $f_\ell$ defines the conditional distribution of $\bY_t$ given $X_t=\ell$. 
The density of $\by$ is defined by
\begin{equation}\label{eq:distfull}
p(\by) = \sum_{\bx \in \{1,\ldots,L\}^{n+1}} \pi_{x_1}  f_{x_1}(\by_1)\prod_{t=1}^n A[x_{t},x_{t+1}] f_{x_{t+1}}(\by_{t+1}).
\end{equation}
 The conditions made on the transition matrix and on the emission distributions are stated by the following set of assumptions. Note that these assumptions are mild and have been considered already in \citet{gassiat2016inference} to state the identifiability of an HMM based on the distribution of three consecutive observations  (see also \citet{de2016minimax,de2017consistent}).

 \begin{assumption}\label{ass:id}
\begin{itemize}
\item The transition matrix $\bA$ has full rank, is irreducible and aperiodic with stationary distribution $\bpi = (\pi_1,\ldots, \pi_L)^\top$. \label{ass:transition}
\item The  densities defining the emission distributions $\{f_{\ell}\}_{\ell=1}^L$ are linearly independent (\emph{i.e.,} if $\boldsymbol{\xi}=(\xi_1,\ldots,\xi_L) \in\mathbb{R}^L$ is such that for any $\bz\in\mathbb{R}^d$, $\sum_{\ell=1}^L \xi_{\ell} f_{\ell}(\bz)=0$ then $\boldsymbol{\xi}=\boldsymbol{0}$) and are square integrable on $\mathbb{R}^d$. \label{ass:emission}
\end{itemize}
\end{assumption}
Under Assumption~\ref{ass:id}, the identifiability of the finite and infinite parameters of a HMM can be obtained from the distribution of three consecutive observations \citep{gassiat2016inference} or from the distribution of a pair of consecutive observations when the emission distributions are defined as translations of the same distribution \citep{gassiat2016nonparametric}. 

 The aim is to make inference on the order $L$. This can be achieved by using the distribution of a pair of consecutive observations. From \eqref{eq:distfull}, the distribution of a pair of consecutive observations $(\bY_t^\top,\bY_{t+1}^\top)^\top$ is defined by the density
\begin{equation}\label{eq:distcouple}
p(\by_t,\by_{t+1}) = \sum_{\ell=1}^L \pi_\ell f_\ell(\by_t) g_{\ell}(\by_{t+1}),
\end{equation}
where $g_\ell$ is the density of $\bY_{t+1}$ given $X_t=\ell$ and is defined by
\begin{equation}\label{eq:cond}
g_\ell(\by_{t+1}) = \sum_{m=1}^L A[\ell,m]f_{m}(\by_{t+1}).
\end{equation}
Note that \eqref{eq:distcouple} is a mixture model where the density of each of the $L$ components is defined as a product of two specific densities. The mixture proportions correspond to the probabilities of latent states defined by the stationary distribution of the Markov chain. Moreover, due to the structure of the HMM, the second density of any component (\emph{i.e.,} $g_\ell$) is a convex combination of the first densities of all the components (\emph{i.e.,} $f_1,\ldots,f_L$), while in the i.i.d. setting $g_\ell$ and $f_\ell$ are not related. Note that the pairs of consecutive observations $(\by_t^\top,\by_{t+1}^\top)^\top$  are identically distributed according to \eqref{eq:distcouple}  but they are not independent due to the dependency between the elements of the whole vector $\by$. The following lemma shows that the order of the HMM can be identified from the distribution of a pair of consecutive observations. 
 \begin{lemma}\label{lemma:ident}
If Assumption~\ref{ass:id} holds true, then $L$ is identifiable from the distribution of a pair of consecutive observations defined by \eqref{eq:distcouple}.
\end{lemma}
As a direct consequence of Lemma~\ref{lemma:ident}, estimating the number of latent states is equivalent to estimating the number of components in \eqref{eq:distcouple}. A specific integral operator can be used to select the number of components in \eqref{eq:distcouple} inspired from \citep{kwon2021estimation}. In this paper, we present some extensions of this approach that permit do deal with the non-independence between the pairs of consecutive observations $(\by_t^\top,\by_{t+1}^\top)^\top$ and to define all the tuning parameters with a new data-driven procedure.



\subsection{Integral operators}
 Let $\spa$ be the Hilbert space of square integrable functions on $\mathbb{R}^d$. We consider the integral operator $T:\spa\to\spa$ defined, for any function $\omega\in\spa$, by
 $$
 [T(\omega)](\bz_{2}) = \int_{\mathbb{R}^d} \omega(\bz_{1}) p(\bz_1,\bz_{2})  \mathrm{d}\bz_{1},
 $$
 where $p$ is the joint distribution given in \eqref{eq:distcouple}.
 From the observed sample $\by$, controlling the accuracy of the estimators of the singular values of $T$ is a delicate task because the density $p(\bz_1,\bz_{2})$ cannot be estimated without bias by the usual kernel method. Therefore, we introduce a smoothed version of the integral operator, denoted by $T_h$, for which we will be able to compute unbiased estimators of its singular values (see Section~\ref{sec:estim}). The operator $T_h:\spa\to\spa$ is defined, for any function $\omega\in\spa$, by
 $$
 [T_h(\omega)](\bz_{2}) = \int_{\mathbb{R}^d} \omega(\bz_{1}) p_h(\bz_1,\bz_{2})  \, \mathrm{d}\bz_{1},
 $$
 where $p_h$ is the function obtained by the convolution between the density $p$ of a pair of consecutive observations given in \eqref{eq:distcouple} and a multivariate kernel defined as a product of univariate kernels, as follows,
 $$
  p_h(\bz_1,\bz_{2})  = \int_{\mathbb{R}^d \times \mathbb{R}^d} p(\by_1,\by_{2})  K_h^d(\bz_1-\by_1) K_h^d(\bz_{2}-\by_2) \, \mathrm{d}\by_1 \mathrm{d}\by_2,
 $$
 where $\bz_1\in\mathbb{R}^d$, $\bz_2\in\mathbb{R}^d$, $K_h^d(\bu)=\prod_{j=1}^d K_{hj}(u_j)$, $\bu=(u_1,\ldots,u_d)^\top \in\mathbb{R}^d$, $K_{hj}$ being univariate kernels and $h>0$ the associated bandwidth.
 Under usual assumptions on the kernel (see Assumption~\ref{ass:kernel}), the ranks of $T$ and $T_h$ are equal to the order of the HMM (see Proposition~\ref{prop:rank}).
\begin{assumption}\label{ass:kernel}
 Each of the kernels $K_{hj}$, for $j=1, \ldots, d$, has a non-vanishing Fourier transform, belongs to $L^1(\mathbb{R}^d)\cap L^2(\mathbb{R}^d)$ and satisfies $ \displaystyle \int u K_{hj}(u) \mathrm{d}u=0$ and $0< \displaystyle \int u^2K_{hj}(u) \mathrm{d}u<\infty$.
\end{assumption}
\begin{proposition}[Proposition~2.1  and Proposition~2.2 in \citet{kwon2021estimation}]\label{prop:rank}
Under  Assumption~\ref{ass:id},
$$
\text{rank}(T)=L,
$$
where rank($T$) is defined as the dimension of the operator $T$. If in addition, Assumptions \ref{ass:kernel} hold true, then
\begin{equation*}
\text{rank}(T_h)=L.
\end{equation*}
\end{proposition}

 Proposition~\ref{prop:rank} implies that the operators $T$ and $T_h$ are compact and admit a singular value decomposition based on $L$ non-zero singular values $\sigma_1(T)\geq\ldots\geq\sigma_L(T)>0$ and $\sigma_1(T_h)\geq\ldots\geq\sigma_L(T_h)>0$, where $\sigma_j(T)$ denotes the $j$-th largest singular value of operator $T$. Hence, for any $j>L$, $\sigma_j(T)=\sigma_j(T_h)=0$. Therefore, estimating the number of latent states $L$ can be achieved by estimating the number of non-zero singular values of $T$. Under regularity conditions on the density of a pair of consecutive observations (see Assumption~\ref{ass:regularity}),   the differences between the non-zero singular values of $T$ and $T_h$ can be controlled (see Lemma~\ref{lemma:controlSVD}).

\begin{assumption} \label{ass:regularity}
 The density function $p$ has partial derivatives at least until order 3 that belong to $L^1(\mathbb{R}^d)\cap L^2(\mathbb{R}^d)$.
\end{assumption}
\begin{lemma}\label{lemma:controlSVD}
Under Assumptions~\ref{ass:id}, \ref{ass:kernel} and \ref{ass:regularity}, we have
$$
\sum_{\ell=1}^L (\sigma_\ell(T) - \sigma_{\ell}(T_h))^2 = O(h^4).
$$
\end{lemma}
Thus, $T$ and $T_h$ have $L$ non-zero singular values,  by Lemma~\ref{prop:rank}, and  one can control the bias induced by the approximation of the non-zero singular values $\sigma_{\ell}(T)$ by the $\sigma_{\ell}(T_h)$, for each $\ell=1,\ldots,L$ and for any bandwidth $h$,  by Lemma~\ref{lemma:controlSVD}. 

\subsection{Estimator of the order of the HMM} \label{sec:estim}
From the observed sample $\by$, we can compute the unbiased estimator of $p_h$  denoted by $\hat{p}_{h,\by}$ defined for any $\bz_1\in\mathbb{R}^d$  and $\bz_2\in\mathbb{R}^d$ by
\begin{equation}\label{eq:phat}
\hat{p}_{h,\by}(\bz_1,\bz_2) = \frac{1}{n} \sum_{t=1}^{n} K_h^d(\bz_1-\by_t)K_h^d(\bz_2-\by_{t+1}).
\end{equation}
Thus, we can deduce the empirical version of the smoothed operator $\hat{T}_{h,\by}$ defined by
$$
\left[\hat{T}_{h,\by} (\omega)\right](\bz_2)=\int \omega(\bz_1) \hat{p}_{h,\by}(\bz_1,\bz_2) \, \mathrm{d}\bz_1.
$$
To estimate the order of the HMM, it suffices to estimate the singular values of $T_h$ by considering a singular value decomposition of $\hat{T}_{h,\by}$.  However, the rank of $\hat{T}_{h,\by}$ is not equal to $L$, since in general, the number of non-zero singular values of such an operator is $n$. Therefore, to build the estimator  $\hat{L}(\tau_{\alpha,h},h)$, we need to apply, on $\sigma_1(\hat{T}_{h,\by}),\ldots,\sigma_{n}(\hat{T}_{h,\by})$, a threshold $\tau_{\alpha,h}>0$ that depends on the probability $\alpha$ of overestimating the order of the HMM (see Section~\ref{sec:th}) and on the bandwidth $h$.  
This estimator of the number of latent states $\hat L$ is defined by
\begin{equation}\label{eq:estimator}
\hat{L}(\tau_{\alpha,h},h) = \text{card}\left(\left\{\ell: r_{\ell}(\hat{T}_{h,\by}) > \tau_{\alpha,h}\right\}\right),
\end{equation}
where for any operator $\mathcal{T}$ we have
\begin{equation}\label{eq:defrj}
r_\ell(\mathcal{T}) = \left[\sum_{j= \ell}^n \sigma_{j}^2(\mathcal{T}) \right]^{1/2}.
\end{equation}
The threshold $\tau_{\alpha,h}$ depends on  the probability  of overestimating the order and the bandwidth. The next section shows that its consistency can be stated   with   suitable choices of $\alpha$ and $h$.


\section{Properties of the estimator of the HMM order} \label{sec:th}

\subsection{Non-asymptotic results}
The following theorem gives an upper-bound on the probability of overestimating the number of latent states when this number is estimated by $\hat{L}(\tau,h)$. This result is stated by controlling $\|\hat{T}_{h,\by} - T_h\|_{HS}=[\int_{\mathbb{R}^{2d}} \left(\hat{p}_h(\bz_1,\bz_2) - p_h(\bz_1,\bz_2)\right)^2d\bz_1d\bz_2]^{1/2}$ (see Lemma~\ref{lemma:eqnorms} in Appendix for the definition) via a concentration inequality. This control is achieved under mild assumptions (Assumptions~\ref{ass:id} and \ref{ass:kernel}) because it only requires that $\text{rank}(T) = \text{rank}(T_h)=L$. The second part of the theorem shows that, under additional conditions, $\hat{L}(\tau_{\alpha,h},h)$ does not underestimate the order of the HMM. Thus, we obtain a lower-bound on the probability that $\hat{L}(\tau,h)=L$.

 \begin{theorem}\label{thm:finitesample}
Under Assumptions~\ref{ass:id} and \ref{ass:kernel}, for any $0<\alpha<1$, there exists some positive threshold  $\tau_{\alpha,h}$ such that the probability to overestimate the number of states is less than $\alpha$, leading that
\begin{equation} \label{eq:overestimation_nknown}
\mathbb{P}(\hat{L}(\tau_{\alpha,h},h)>L) < \alpha,
\end{equation}
with 
\begin{equation} \label{eq:tau}
\tau_{\alpha,h}= \frac{\|K_h\|_2^{2d}}{n^{1/2}} \left[\left(\frac{n+1}{n}C_{\alpha,1}\right)^{1/2}  + C_{2} ^{1/2}\right],\end{equation}
where $\|K_h\|_2^2=\int_{\mathbb{R}} K_h^2(u)du$, $C_{\alpha,1}=36 \ln(1/\alpha) t_{\text{mix}} $, $C_2= 9+8 t_{\text{mix}}$ and $t_{\text{mix}}$ is the mixing time of the underlying Markov chain recalled in Definition \ref{def:tmix}  in the Appendix. 

If in addition Assumption~\ref{ass:regularity} holds true and if $h$ is small enough and $n$ is large enough to ensure that for some $\varepsilon >0$, $\sigma_L(T)>2\tau_{\alpha,h} + \varepsilon$, then
\begin{equation}\label{eq:underestimation_nknown}
\mathbb{P}(\hat{L}(\tau_{\alpha,h},h)<L) =0
\text{ and }
\mathbb{P}(\hat{L}(\tau_{\alpha,h},h)=L) \geq 1-\alpha.
\end{equation}
\end{theorem}
From \eqref{eq:overestimation_nknown}, the probability of overestimating the order of the HMM can be set as small as wanted for any value of the bandwidth $h$, by considering the threshold given by \eqref{eq:tau}. Therefore, even if the singular values of $T$ are estimated with bias from $T_h$ (when $h>0$, $\text{rank}(T) = \text{rank}(T_h)=L$ but for any $j=1,\ldots,L$ $\sigma_j(T)\neq\sigma_j(T_h)$), we can make the probability of overestimating $L$ as small as wanted. However, to avoid underestimating $L$, the method requires to consistently estimate $p$ by $p_h$, and so that the bandwidth $h$ tends to zero at a suitable rate in order to have $\sigma_L(T)>2\tau_{\alpha,h}$.  Thus, only the variance of the estimators $\sigma_j(\hat{T}_{h,\by})$ can lead to the overestimation of $L$. Indeed, despite the bias, the ranks of $T$ and $T_h$ are the same. However, both the bias and the variance of the estimators $\sigma_j(\hat{T}_{h,\by})$ impact the underestimation. Note that the condition $\sigma_L(T)>2\tau_{\alpha,h}$ cannot be verified in practice since it depends on  the singular values of the theoretical operator. The following section gives  rules on $\alpha$ and $h$, which are sufficient to ensure the consistency of the estimator \eqref{eq:estimator}.

We now  discuss the connections of the results stated by Theorem~\ref{thm:finitesample} and those presented in Theorem 3.1 in  \citet{kwon2021estimation}. Both theorems allow for a control of the probability of overestimating the rank of the operator by controlling the concentration of $\|\hat{T}_{h,\by} - T_h\|_{HS}$. However, this control is achieved in two different manners. In the i.i.d. setting, see Proposition 3.1 of their paper, a concentration inequality for $\|\hat{T}_{h,\by} - T_h\|_{HS}$ is obtained by combining Theorem 3.4 of \citet{pinelis1994optimum} applied to sums of independent random elements in the space of Hilbert-Schmidt operators and Hoeffding's concentration inequality. In the latter, the concentration bound, and thus the threshold, can be explicitly bounded as it involves quantities which only depend on the kernel and the bandwidth, see (3.8) in their paper.
Similar reasoning cannot be used in our context due to the dependency between the observations implied by the HMM structrure. Thus, the proof of Theorem~\ref{thm:finitesample} presented in Appendix relies on specific statistical tools for HMM since it combines McDiarmid inequalities for HMM and coupling methods. The bound obtained in Theorem~\ref{thm:finitesample} involves the mixing time of the Markov chain, which is unknown and cannot be estimated without any knowledge on the order of the HMM. Thus, contrary to the i.i.d. context, we do not bound the quantities involved in the concentration inequality. This choice raises the question of the tuning of some constant that we discuss in Section~\ref{sec:cst}.

\subsection{Asymptotic results}
The following corollary states the consistency of the estimator of the number of states defined by \eqref{eq:estimator}. This consistency is obtained by considering an appropriate rate of decreasing of the probability $\alpha_n$ of overestimating the order and a suitable bandwidth $h_n$ whose values depend on the sample size $n$, without requiring  the consistency on the smallest non-zero singular value of $T$. 
 The threshold $\tau_{\alpha_n,h_n}$ depends on both the probability of overestimating the order and the bandwidth (see \eqref{eq:tau}). The quantities $\alpha_n$ and $h_n$ tend to zero when $n$ tends to infinity. 
 However, these quantities should tend to zero at an appropriate rate which ensures that $\lim_{n\to\infty}$ $r_{L+1}(\hat{T}_{h_n,\by})/\tau_{\alpha_n,h_n}=0$ and $\lim_{n\to\infty}$ $\tau_{\alpha_n,h_n}/r_{L}(\hat{T}_{h_n,\by})=0$. 
 Indeed, considering $\alpha_n$ tending to zero ensures that the order is not overestimated. However, to ensure that the order is not underestimated (see condition $\sigma_L(T)>2\tau_\alpha + \varepsilon$ in Theorem~\ref{thm:finitesample}), the threshold needs to tend to zero as the sample size tends to infinity.

\begin{corollary}\label{thm:asymptotic}
Under Assumptions~\ref{ass:id} and \ref{ass:kernel},  if it exists $u>0$ such that $\|K_h\|_2^2 \lesssim h^{-u}$, then for a bandwidth $h_n=O(n^{-\beta})$ with $0<\beta<(2du)^{-1}$, then considering the threshold $\tau_{\alpha_n, h_n}$
satisfying
\begin{equation} \label{eq:thresholdv}
\tau_{\alpha_n, h_n} = o(1) \text{ and }  n^{-\frac{1}{2} + du\beta} \tau_{\alpha_n, h_n}^{-1}=o(1),
\end{equation}
 implies that $\hat{L}(\tau_{\alpha_n, h_n},h_n)$ is a consistent estimator of $L$ meaning that
\begin{equation*} 
\lim_{n\to \infty} \mathbb{P}(\hat{L}(\tau_{\alpha_n,h_n},h_n)=L) =1.
\end{equation*}
\end{corollary}
The conditions \eqref{eq:thresholdv} in Corollary~\ref{thm:asymptotic} permit to avoid underestimation and overestimation of the order of the HMM. Indeed, since $\tau_{\alpha_n, h_n} = o(1)$, for $n$ large enough $\sigma_L(T)>2\tau_{\alpha_n,h_n} + \varepsilon$ and thus the approach does not asymptotically underestimate the order of the HMM (see \eqref{eq:underestimation_nknown} in Theorem~\ref{thm:finitesample}). Moreover, 
since $n^{-\frac{1}{2} + du\beta} \tau_{\alpha_n, h_n}^{-1}=o(1)$, then $\left(\frac{n+1}{n}C_{\alpha_n,1}\right)^{1/2}  + C_{2} ^{1/2}$ tends to infinity as $n$ tends to infinity (see \eqref{eq:tau} in Theorem~\ref{thm:finitesample}), and so $\alpha_n$ tends to zero leading that the  approach does not asymptotically overestimate the order of the HMM. Corollary~\ref{thm:asymptotic} implies the consistency of the estimator by defining the threshold
\begin{equation} \label{eq:threshold2}
\tau_{\alpha_n, h_n} = C n^{-\frac{1}{2} + du\beta} \ln n,
\end{equation}
for any positive constant $C$. If the choice of the value of the constant does not influence the asymptotic behavior $\hat{L}(\tau_{\alpha_n,h_n},h_n)$ because this threshold respects \eqref{eq:thresholdv}, the impact of the unknown constant $C$ on the resulting estimator can be strong on finite sample size. Hence, we now discuss computational aspects of the method, including the tuning of this constant.

\section{A Data-driven calibration procedure}\label{sec:cst}
In this section, we present a data driven procedure  to calibrate  positive constant $C$ in a threshold procedure based on  the threshold defined by:
$$
\tau_n=C \upsilon(n),
$$
where $\upsilon(n)$ is a positive function that tends to zero when $n$ tends to infinity. This situation covers the problem of the proposed procedure and its threshold \eqref{eq:threshold2} by considering $\upsilon (n)=n^{-\frac{1}{2} + du\beta} \ln n$.
The idea of the procedure is to split the data into $S$ mini-batches of size $m$ and to compute the statistics requiring thresholding on each mini-batch. Then, a K-means is run with two groups on the means of these statistics computed over the $S$ mini-batches. We show that, under some assumptions made on the first $L$ statistics, the K-means approach allows to split the significant and the non-significant statistics. Thus, we select the constant such that the thresholding rule applied on the mini-batches provides the same estimator. Then, we use the thresholding procedure on the whole sample by considering this particular value of the constant. The proposed method can be run to tune the constant in  \eqref{eq:threshold2} but also in the spectral-approach proposed by \citep{lehericy2019consistent}. Numerical experiments presented in Section~\ref{sec:ex} illustrate the interest of the proposed approach.

Considering the  mini-batch $\by^{[s]}$ of size $m$, we compute the statistics $r_1(\widehat{T}_{m^{-\beta},\by^{[s]}}),\ldots,r_{m}(\widehat{T}_{m^{-\beta},\by^{[s]}})$ and we define for $\ell=1,\ldots,m$
$$
\widehat{r}_\ell = \frac{1}{S}\sum_{s=1}^Sr_\ell(\widehat{T}_{m^{-\beta},\by^{[s]}}),
$$
$S$ being the number of mini-batches of size $m$. We perform a clustering of the values $\widehat{r}_1,\ldots,\widehat{r}_{m}$ into two groups in order to minimize the within-group variance with K-means. The idea behind this clustering is to group all the values of $\widehat{r}_\ell$ with $\ell>L$ into the same cluster since these values should concentrate around zero. Let $g_1$ be the cluster with the highest mean and $G_1$ the subset of $\{1,\ldots,m\}$ belonging to $g_1$. The following proposition gives sufficient conditions that ensures that the K-means algorithm run in two groups on the  $\widehat{r}_1,\ldots,\widehat{r}_m$ groups the first $L$ values into the same groups. Therefore, under assumptions on $\widehat{r}_1,\ldots,\widehat{r}_L$, if the size of the mini-batches is large enough then $\wwr_{\ell}$ is as small as wanted with high probability when $\ell>L$ and so the cardinal of $g_1$ is equal to the $L$.
\begin{proposition} \label{prop:Kmeans}
Let $\wwr_1,\ldots,\wwr_m$ be positive variables with $\wwr_{\ell}\geq \wwr_{\ell+1}$ and $\wwr_{L+1}<\varepsilon$ for some positive $\varepsilon$. Then, if $\varepsilon$ is small enough, $m$ is large enough and if $(1/s) \sum_{\ell=1}^s \wwr_\ell<[(s/\xi(s+1))^{1/2} + 1] \wwr_{s+1}$ for any $s\in\{1,\ldots,L-1\}$, where $\xi>1$, then K-means run with two groups gathers all the  first $L$ variables in the same group, if $m$ is fixed.
\end{proposition}

Considering the partition provided by the K-means, a threshold that allows the groups provided by the clustering to be recovered is
$$
\sigma_m := \max_{\ell \in \{1,\ldots,m-1\}\setminus G_1} \widehat{r}_\ell .
$$
The idea behind the definition of $\sigma_m$ is to set the smallest threshold that would provide a number of latent states equal to the cardinal of $g_1$ if the method would be applied on the values of $\widehat{r}_1,\ldots,\wwr_m$. To set the constant, it suffices to consider that the threshold $\sigma_m$ has the shape $\widehat{C} \upsilon(m)$, leading that
\begin{equation*}
   \widehat{C} =  \sigma_m / \upsilon(m).
\end{equation*}
Therefore, when $\widehat{C}$ has been tuned on the mini-batches, the estimation of the order can be performed on the full sample by  plugging this constant  into the threshold \eqref{eq:thresholdGauss}.

As a consequence of Proposition~\ref{prop:Kmeans}, the strategy used to tune the constant is relevant since $\wwr_{L+1}$ converges to zero as $m$ tends to infinity. Indeed, this strategy allows for a consistent detection of the order of the HMM under the assumptions made on the first $L$ statistics. Note that because the method is used with fixed sample size of the mini-batches (\emph{i.e.,} $m$ does not grow with $n$), if the assumption made on the significant statistics $\hat{r}_{\ell}$ is not significant then the estimator of $L$ obtained by considering the provided constant is still consistent, despite the fact that the K-means procedure is not consistent to estimate $L$. Indeed, the role of the K-means procedure is only to provided a relevant value for the constant in \eqref{eq:threshold2}. 
With a careful reading of the proof, we can see that the K-means procedure is still consistent if the size of the mini-batches $m$ grows with the sample size $n$ such that $m/\upsilon(n)$ tends to zero as $n$ tends to infinity. In such case, consistency of the K-means procedure only requires that the assumption made on the first $L$ statistics is satisfied. As an example, considering continuous data and a Gaussian kernel leads that $u=1$ (since $\|K_h\|_2^2=O(h_n^{-1})$) where we consider the same bandwidth $h=n^{-\beta}$ for the kernel density estimation with the usual bandwidth $\beta=1/(4+2d)$ after scaling the data. Thus, in such case, the proposed threshold is defined by 
\begin{equation}\label{eq:thresholdGauss}
    \hat{\tau}_{\alpha_n,h_n}= \widehat{C}  n^{-\frac{1}{2+d}} \ln n,
\end{equation}
where $\widehat{C}$ is defined by the data-driven procedure based on K-means.

\section{Computing the singular values} \label{sec:computation}

To estimate the rank of $T_h$, it suffices to estimate the singular values of $T_h$ by considering a singular value decomposition of $\hat{T}_{h,\by}$. However, performing the singular value decomposition (SVD) of an operator directly is not straightforward. Therefore, we introduce a  $n\times n$ matrix $\widehat{\bV}_{h,\by}$ that has the same singular values as $\hat{T}_{h,\by}$ and for which we will be able to perform the SVD.
Let the empirical $n\times n$ matrix $\widehat{\bV}_{h,\by}$ be defined by
\begin{equation}\label{eq:matV}
\widehat{\bV}_{h,\by}=\frac{1}{n}\Delta_1^\top \widehat{\bW}_{h,\by}^{1/2}\Delta_1\Delta_{n+1}^\top \widehat{\bW}_{h,\by}^{1/2}\Delta_{n+1},
\end{equation}
where 
$\Delta_1$ and $\Delta_{n+1}$ are the $(n+1)\times n$ matrices defined as block matrices by
$
\Delta_1^\top = \begin{bmatrix}
\boldsymbol{0}_{n}
\, \boldsymbol{I}_{n} 
\end{bmatrix}$
and
$\Delta_{n+1}^\top = \begin{bmatrix}
\boldsymbol{I}_{n} 
\, \boldsymbol{0}_{n}
\end{bmatrix}
$, where $\boldsymbol{I}_{n}$ is the identity matrix of size $n\times n$ and $\boldsymbol{0}_{n}$ is the null vector of length $n$, and where $\widehat{\bW}_{h,\by}$ is the $(n+1)\times (n+1)$ matrix defined by
$
\widehat{\bW}_{h,\by}[t,s] = \phi_{h}(\by_t,\by_{s}),
$
for $1\leq t,s\leq (n+1)$ where the function $\phi_h$ is such that $\phi_h: \mathbb{R}^{2d} \to \mathbb{R}$ and
$$
\phi_{h}(\by_t,\by_{s})=\int_{\mathbb{R}^d} K_h^d(\bz_1 - \by_t) K_h^d(\bz_1-\by_{s}) \, \mathrm{d}\bz_1,
$$
where $h>0$ is a bandwidth and $K^d$ is a $d$-dimensional kernel. In our setting, the kernel $K_h^d$ is defined as a product of univariate kernels, $K_h^d(\bu)=\prod_{j=1}^d K_{hj}(u_j)$, which permits to rewrite 
$\phi_{h}(\by_t,\by_{s})= \prod_{j=1}^d \int_{\mathbb{R}} K_{hj}(z_{1j} - y_{tj}) K_{hj}(z_{1j}- y_{sj}) \, \mathrm{d}z_{1j}. $
In many cases, the function $\phi_h$ can be computed in closed form.  For instance, with the Gaussian kernel defined on $\mathbb{R}$ by $K_{hj}(x)= (2h^2\pi)^{-1/2} \exp ( -x^2/2h^2)$, we obtain $\phi_h(\by_t,\by_{s})= (2h\sqrt{\pi})^{-d} \exp (- \sum_{j=1}^d (y_{tj}- y_{sj})^2 / (4h^2))$, and with the Von-Mises Kernel defined on $[-\pi, \pi]$ by $K_{hj}(x)= (2 \pi I_0(h^{-2}))^{-1} \exp ( \cos(x)/h^2)$, with $I_0$ being the modified Bessel function of the first kind, we obtain that 
$\phi_h(\by_t,\by_{s})= (\sqrt{2 \pi} I_0(h^{-2}))^{-2d} \prod_{j=1}^d I_0(2 \cos((y_{tj}- y_{sj})/2) / h^2)$. 

As a direct consequence of Corollary~3.1 of \citet{kwon2021estimation}, the singular values of $\bV_{h,\by}$ are equal to the singular values of $\hat{T}_{h,\by}$ leading that for any $j=1,\ldots,n$, we have
$$
\sigma_j(\hat\bV_{h,\by}) = \sigma_j(\hat{T}_{h,\by}).
$$


\section{Numerical experiments} \label{sec:ex}
During the numerical experiments, the results of proposed method are obtained by considering the threshold \eqref{eq:thresholdGauss} and have been obtained by the R package \pkg\ used with the computational aspects described in Section~\ref{sec:computation} with $m=125$, and with Gaussian kernel for the continuous data and Von-Mises kernel for the circular data. Data are generated from an HMM defined with $L=3$ hidden states and transition matrix
$$
A_{\nu} = \begin{bmatrix}
1 - 2\nu & \nu & \nu \\
\nu & 1-2\nu & \nu \\
\nu & \nu & 1-2\nu
\end{bmatrix},
$$
where the parameter $\nu$ allows us to define different mixing times. Indeed, in the case where $\nu=1/3$, this setup generates independent data, while the mixing time increases when $\nu$ tends to zero. Conditionally on the hidden state $X_t$, the components of the vector $\bY_t=(Y_{t1},\ldots,Y_{td})$ are independently generated from the model defined, for any $(t,j)$, by $
Y_{tj} =  \left(\mathds{1}_{\{X_t=2\}} - \mathds{1}_{\{X_t=3\}}\right)\delta + \varepsilon_{tj},
$
where all $\varepsilon_{tj}$ are generated independently and $\delta$ is a constant tuning parameter. Three distributions for continuous data are considered for $\varepsilon_{tj}$ (Gaussian, Student with three degrees of freedom, Laplace) and one distribution for circular data is considered for $\varepsilon_{tj}$  (Von-Mises). The parameter $\delta$ allows us to tune the overlaps between the emission distributions of each state.

In this experiment, we compare our integral-based method and the spectral method for recovering the true order $L=3$. The spectral method estimates the order from the rank of the $(M\times M)$ matrix $\widehat{N}_M$ defined by $\widehat{N}_M[k,\ell]=\frac{1}{n}\sum_{t=1}^n \phi_k(\bY_t)\phi_{\ell}(\bY_{t+1})$ where the functions $\phi_1,\ldots,\phi_M$ are basis functions. The spectral method is used with histogram basis defined by the $M$ quantiles. Moreover, a threshold must be applied on the observed singular values of $\widehat{N}_M$. As suggested in Section 5.3.2 of \citet{lehericy2019consistent}, we consider an integer $M_{\text{reg}}\leq M$ and we estimate the affine dependency of the singular values of $\widehat{N}$ with respect to their index with a linear regression using its $M_{\text{reg}}$ smallest singular values. Then, we set a thresholding parameter $\tau=1.5$ and we say that a singular value is significant if it is above $\tau$ times the value that the regression predicts for it. The estimator of the order is the number of consecutive significant singular values starting from the largest one. Thus, the spectral method has three tuning parameters: the basis family, the number of basis elements $M$ (note that by construction, the estimator of the order is upper-bounded by $M$) and the number $M_\text{reg}$ of singular values used for determining the threshold. In our experiments, we use $M=20$ and $Mreg=10$. Alternatively, we use the spectral method where the threshold is defined by $C \sqrt{\ln n/n}$, as suggested in \citet{lehericy2019consistent}. This constant is tune with the proposed approach based on mini-batches and Kmeans algorithm described in Section~\ref{sec:cst}.

The methods are compared by considering different overlaps between the emission distributions. The data are thus generated with different values of $\delta$, which define different error rates (2.5\%, 5\% and 10\%) based on the distribution of a single observation. Table~\ref{tab:simudelta} indicates the order estimated by the proposed method on 100 samples generated with $d=1$,  $\nu=0.1$, different sample sizes $n$, different families of emission distributions and  an  error rate of 5\%. Table~\ref{tab:simudelta0025} and Table~\ref{tab:simudelta01}, presented in Appendix, show the results obtained with error rates of 2.5\% and 10\%. Overall, the results illustrate the consistency of the proposed estimator, as stated by Theorem~\ref{thm:asymptotic}. Indeed, for any family of distributions and any value of $\delta$, the estimator selects the true order with probability $1$ when the sample size increases. Moreover, for  small samples, the estimator does not overestimate the order but can underestimate it. Note that this phenomenon was already observed in the i.i.d. setting. Finally, the more different the emission distributions, the more accurate the estimator for small samples. Furthermore, the comparison made with the spectral method proposed in \citet{lehericy2019consistent} shows that our approach gives the best results whatever the setting and the data-driven procedure proposed in this paper improves his results obtained with the slope heuristics.

\begin{table}[htp!]
\centering
\begin{scriptsize}
\begin{tabular}{cc|ccc c|ccc c|ccc c|ccc c}
  \hline
Method & $n$ & \multicolumn{4}{c|}{Gaussian} & \multicolumn{4}{c|}{Student} & \multicolumn{4}{c|}{Laplace} & \multicolumn{4}{c}{Von-Mises}\\
& & L-1 & L-2 & \normalfont{L-3} & L$>$3 & L-1 & L-2 & \normalfont{L-3} & L$>$3 & L-1 & L-2 & \normalfont{L-3} & L$>$3 & L-1 & L-2 & \normalfont{L-3} & L$>$3 \\ 
  \hline
proposed & 250 & 0 & 40 & 60 & 0 & 0 & 37 & 63 & 0 & 0 & 19 & 81 & 0 & 0 & 7 & 89 & 4 \\ 
    & 500 & 0 & 23 & 77 & 0 & 0 & 15 & 85 & 0 & 0 & 6 & 94 & 0 & 0 & 1 & 97 & 2 \\ 
    & 1000 & 0 & 0 & 100 & 0 & 0 & 0 & 100 & 0 & 0 & 0 & 100 & 0 & 0 & 0 & 100 & 0 \\ 
    & 2000 & 0 & 0 & 100 & 0 & 0 & 0 & 100 & 0 & 0 & 0 & 100 & 0 & 0 & 0 & 100 & 0 \\ 
    & 4000 & 0 & 0 & 100 & 0 & 0 & 0 & 100 & 0 & 0 & 0 & 100 & 0 & 0 & 0 & 100 & 0 \\ 
   \hline
  spectral & 250 & 0 & 0 & 4 & 96 & 0 & 0 & 1 & 99 & 0 & 0 & 0 & 100 & 0 & 1 & 2 & 97 \\ 
slope & 500 & 0 & 1 & 11 & 88 & 0 & 0 & 1 & 99 & 0 & 0 & 0 & 100 & 0 & 0 & 4 & 96 \\ 
  & 1000 & 0 & 0 & 7 & 93 & 0 & 0 & 0 & 100 & 0 & 0 & 0 & 100 & 0 & 0 & 5 & 95 \\ 
   & 2000 & 0 & 0 & 11 & 89 & 0 & 0 & 0 & 100 & 0 & 0 & 0 & 100 & 0 & 0 & 7 & 93 \\ 
   & 4000 & 0 & 0 & 15 & 85 & 0 & 0 & 0 & 100 & 0 & 0 & 0 & 100 & 0 & 0 & 3 & 97 \\ 
   \hline
  spectral & 250 & 18 & 11 & 24 & 44 & 13 & 24 & 22 & 38 & 16 & 22 & 27 & 34 & 11 & 23 & 24 & 39 \\ 
kmeans   & 500 & 17 & 26 & 47 & 10 & 11 & 31 & 40 & 18 & 17 & 30 & 37 & 14 & 8 & 42 & 38 & 12 \\ 
  & 1000 & 4 & 18 & 75 & 3 & 4 & 22 & 66 & 8 & 7 & 19 & 69 & 5 & 4 & 19 & 70 & 7 \\ 
  & 2000 & 0 & 0 & 100 & 0 & 0 & 2 & 98 & 0 & 0 & 4 & 96 & 0 & 1 & 2 & 97 & 0 \\ 
   & 4000 & 0 & 0 & 100 & 0 & 0 & 0 & 100 & 0 & 0 & 0 & 100 & 0 & 0 & 0 & 100 & 0 \\ 
   
     \hline
\end{tabular}
\end{scriptsize}
\caption{Percentage of number of states selected by the competing methods (proposed method "proposed", spectral method with slope heuristic used for tuning threshold "spectral slope" and spectral method with the proposed method based on Kmeans used for tuning the constant "spectral kmeans"), according to the family of the emission distributions and the sample size, obtained on 100 replicates generated with $d=1$, $\nu=0.1$, with an marginal overlap between the emission distributions of 5\%. }\label{tab:simudelta}
\end{table}

\section{Identifying movement regimes from masked boobies trajectories} \label{sec:appli}

Following  \cite{nathan2008movement},  animal movement depends on internal states of individuals and therefore the segmentation of the trajectories helps ecologists to characterize different movement patterns interpreted as different internal states. We present one example of such an approach by studying  three trajectories of a masked booby bird living on Meion Island in the Fernando de Noronha archipelago in Brazil. The bird is equipped with a GPS that records its position every 10 seconds. The initial data expressed in latitude and longitude are projected onto a local coordinate system and we use  a HMM approach to identify different movement patterns within the trajectory. The number of patterns is unknown and \cite{pohle2017selecting} advocate that AIC and BIC tend to overestimate the order of the HMM and that this order should be chosen according to biological consideration. We compare the order estimated by the different criteria and the order estimated using the method proposed in this paper.   

The data consist in three different trips composed of 2712, 2451, and 2229 GPS positions respectively. The data come from the field work of Sophie Bertrand (IRD), Guilherme Tavares (UFRGS), Christophe Barbaud and Karine Delord (CNRS) and are kindly made available by the IRD Tabasco Young International Team (JEAI).

Rather than analyzing absolute positions, the classical movement ecology approach consists in deriving different metrics from the sequence of positions. This either the bivariate sequence $(L_t, \phi_t)_{ 1 \leq t\leq n}$ of the step length  and turning angle sequence as proposed in \cite{vermard2010identifying, walker2010pioneer}  or the bivariate sequence $(V^p_t,V^r_t)_{ 1 \leq t\leq n}$ of the orthogonal components of persistence velocity and rotational velocity  as in \cite{gurarie2009novel}. Biological knowledge concerning the movement of these birds supports the use of these metrics to distinguish between different behaviors.

We used the R package moveHMM \citep{michelot2016movehmm}, a very popular package in the movement ecology community, to fit a HMM on the bivariate  sequence $(L_t, \phi_t)_{ 1 \leq t\leq n}$ where  the emission distributions are assumed to be a product between a gamma distribution and a Von Mises distribution to respect the typical choice of such analysis as in \cite{morales2004extracting}, this model will be referred to as GVM model.  We also fit a Normal Mixture model, referred to as NM model, using depmixS4 \citep{visser2010depmixs4} on the bivariate sequence of persistence velocity and rotational velocity. Table \ref{tab:resRedfooted} presents the order of the HMM obtained by the proposed method and by the information criteria (AIC, BIC and ICL) under the parametric assumptions detailed previously.


\begin{table}[ht!p]
    \centering
    \begin{tabular}{r|cccc}
\hline
Variables & Method & Trip 1 & Trip 2 & Trip 3   \\
\hline
 step length and turning angle & Proposed method & 3 & 3 &  3 \\
& AIC (GVM) & 10 & 9 &  10 \\
& BIC (GVM)& 10 & 5 &  8  \\
& ICL (GVM) & 4 & 6 &  4\\
\hline
velocity and rotational velocity  & Proposed method & 3 & 3 & 3  \\
& AIC (NM) & 10 & 10 & 9  \\
& BIC (NM) & 9 & 7 &  7 \\
& ICL (NM) & 9 & 7 &  7 \\
\hline
\end{tabular}
    \caption{Estimated order of the HMM obtained by the proposed method, by AIC, BIC and ICL for the two parametric models and for three different trips of the same Red-footed booby ({\it Sula sula}) individual.}
    \label{tab:resRedfooted}
\end{table}
We propose to compare the estimators of the latent sequences obtained by considering the variables of step length and turning angle and by considering the variables of  velocity and rotational velocity. Thus, we consider the order selected by the proposed method and when kernel density estimator are used to fit the emission distributions.  
The adjusted Rand indexes between the two estimators of the latent sequences are 0.959, 0.751 and 0.980 for the trips 1, 2 and 3 respectively.

\begin{figure}
    \centering
    \includegraphics[scale=0.4]{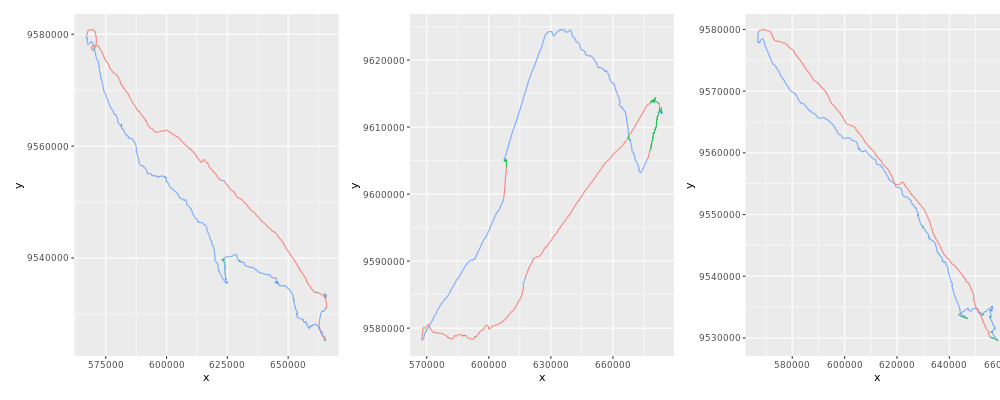}
    \caption{Latent states during the trips obtained by considering the variables of step length and turning angle}
    \label{fig:my_label}
\end{figure}

\begin{figure}
    \centering
    \includegraphics[scale=0.4]{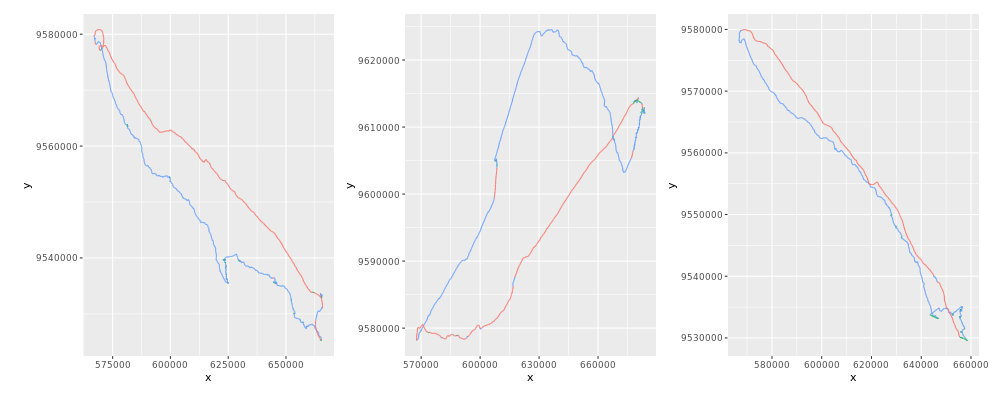}
    \caption{Latent states during the trips obtained by considering the variables of velocity and rotational velocity}
    \label{fig:my_label2}
\end{figure}

It is difficult to justify the choice of one or other family for the emission distributions, however a poorly-adapted choice will often lead to an overestimation of the HMM order.  This is likely the case when fitting the model to the velocity sequence, which explains the very high order estimated by the penalized likelihood criteria. In contrast, the non-parametric approach finds the same order for each of the bird's trips considered and for the two metrics derived from the initial GPS relocation data. On this example, the non parametric approach proposed in this paper seems more robust and provides an interesting alternative to parametric views, which tends to favor high orders, as indicated in \citet{pohle2017selecting}. For the two metrics considered, our approach returns three states (see Figure \ref{fig:my_label} and Figure \ref{fig:my_label2}), which corroborates the interpretation of the ecologists. State $1$ corresponds to the activity with the lowest average speed, state $2$ with a medium speed, while state $3$ corresponds to the fastest average speed.

\section{Conclusion} \label{sec:cl}
In this paper, we introduced a new estimator for assessing the order of a non-parametric HMM by using the rank of an integral operator. To take into account the variability of the empirical singular values of a smoothed version of this operator, a data-driven thresholding rule is proposed and based on a new heuristics for setting the unknown constants. We give a theoretical justification and numerical studies have shown that it can be successfully applied to several problems.
Under standard assumptions for non-parametric HMMs (\emph{i.e.,} full rank covariance matrix and linear independence between the emission distributions), consistency of the estimator is established. 
As illustrated on benchmark and real data, the proposed approach considers different types of data including continuous data but also multivariate data or directional data. Numerical experiments illustrate that the proposed approach gives good results and provides greater flexibility enabling the modeling and analysis of more complex data. 



\bibliography{biblio.bib}

\appendix

\section{Preliminary results}\label{Preliminary:results}

We recall here some definitions and basic properties of the mixing rate of Markov chains. We also provide some useful Lemmas, which will be proven in Appendix~\ref{App:proofl}. 

In the following, we denote by $(X_t)_{t\in \mathbb{Z}}\sim (A,\bpi)$ a Markov chain with irreducible transition matrix $A$ and stationary distribution $\bpi$ and denote by $\vec{X}_t$ the pair $(X_t, X_{t+1})$.  First of all, we define the total variation distance between two probability distributions $\mu$ and $\nu$ on $ \{1,\ldots, L\}$ by

$$\|\mu-\nu \|_{TV}=\frac{1}{2}\sum_{\ell \in \{1,\ldots, L\}}\vert \mu(\ell)-\nu(\ell)\vert.$$

\begin{definition}\label{def:uniform_ergo}
Let $(X_t)_{t\in \mathbb{Z}}\sim (A,\bpi)$ a Markov chain and define, for any $t$, 
\begin{equation*}
d(t)=\sup_{\ell \in \{1,\ldots,L\}} \| A^t(\ell,.)-\bpi \|_{TV}.
\end{equation*}
The Markov chain is \emph{uniformly ergodic} if $d(t)$ goes to zero at a geometric rate as $t$ goes to infinity, \emph{i.e.}
if there exists two constants $c>0$ and  $0< \rho <1$ such that  $d(t) \leq c \rho^t$. 
\end{definition}

However, the constants supplied by Definition \ref{def:uniform_ergo} are generally very difficult to calculate and too conservative to be of any practical use. Fortunately, we can do much better with the spectral gap and the mixing time. Hence, it is useful to introduce the mixing time variable which measures the time required by a Markov chain for the distance to stationarity to be small.

\begin{definition}\label{def:tmix}
Let $(X_t)_{t\in \mathbb{Z}}\sim (A,\bpi)$ . Its mixing time $t_{\text{mix}}$ is defined by
$$t_{\text{mix}}=\min\{t : d(t) \leq \frac{1}{4}\}.$$
\end{definition}

It is well known that irreducible and aperiodic finite state chains are always uniformly ergodic 
and we can easily note that in this case, $t_{\text{mix}}$ is finite and can be bounded by the (pseudo)-spectral gap of the transition matrix, see \cite{wolfer2019estimating}. \\

First of all, we need to ensure that the hidden process $\vec{\bZ}_{s}=(\vec{X}_s, \vec{\bY}_{s})$ is a Markov chain with a mixing time that we are able to control, using the mixing time of the hidden Markov chain $(X_t)_t$.

\begin{lemma}\label{lemma: Mixing_Hidden} 
Suppose that Assumption~\ref{ass:id}  holds true. Let $X_1, \ldots, X_n \sim (A, \bpi)$  be a Markov chain with mixing time $t_{mix}$ and stationary distribution $\bpi$.
The hidden chain $\vec{\bZ}=(\vec{\bZ}_{1}, \ldots, \vec{\bZ}_{n})$  is also then a Markov chain  with kernel transition $A^{\vec{\bZ}}: \{1,\ldots, L\}^2\times \mathbb{R}^{2d} \to [0,1]$ and stationary distribution $\bpi^{\vec{\bZ}}=\bpi^{\vec{\bX}}\otimes G$ with $G$ the transition kernel from $ \{1,\ldots, L\}^2$ to $\mathbb{R}^{2d}$ and $\bpi^{\vec{\bX}}$ such that  $\bpi^{\vec{\bX}}(\ell, k)=\mathbb{P}_{\bpi}(X_{s}=\ell, X_{s+1}=k)$ for all $(\ell,k) \in \{1,\ldots, L\}^2$.  
Furthermore, the mixing time of the Markov chain $\vec{\bZ}$ denoted $t_{mix}^{\vec{\bZ}}$ is at most $t_{mix}+1$.

\end{lemma}

\begin{lemma}[Covariance inequality for hidden Markov chains]\label{lemma: Cov_Paulin_hidden} 
Suppose that Assumptions~\ref{ass:id} and \ref{ass:kernel} hold true. Let $X_1, \ldots, X_n \sim (A, \bpi)$ be a Markov chain with mixing time $t_{mix}$ and stationary distribution $\bpi$ and define the hidden Markov chain $(\vec{\bZ}_{1}, \ldots, \vec{\bZ}_{n}) \sim (A^{\vec{\bZ}},\bpi^{\vec{\bZ}})$ from Lemma \ref{lemma: Mixing_Hidden} with mixing time $t_{mix}^{\vec{\bZ}}$. Then for any measurable function in $L^2(\bpi^{\vec{\bZ}})$: $\phi: \{1,\ldots, L\}^2\times \mathbb{R}^{2d} \to \mathbb{R}$, we have

$$\sum_{t=1}^{n}\mathbb{E}[\phi(\vec{\bZ}_1)\phi(\vec{\bZ}_t)]\leq 4 t_{mix}^{\vec{\bZ}} \mathbb{V}[\phi(\vec{\bZ_1})].$$
In particular,  for $\psi:  \mathbb{R}^{2d} \to \mathbb{R}$ we have

$$\sum_{t=1}^{n}\mathbb{E}[\psi(\vec{\bY}_1)\psi(\vec{\bY}_t)]\leq 4 t_{mix}^{\vec{\bZ}} \mathbb{V}[\psi(\vec{\bY}_1)].$$
\end{lemma}

\begin{lemma}[Hoffman-Wielandt inequality \citep{bha1994}] \label{lemma:Hoffman}
Let  $\mathcal{T}$ and $\mathcal{T}'$ be two operators with finite ranks, then we have that for any positive integer $j$
$$
\sum_{i \geq j} \left( \sigma_i(\mathcal{T}) - \sigma_i(\mathcal{T}')\right)^2 \leq \| \mathcal{T} - \mathcal{T}'\|_{HS}^2.
$$
\end{lemma}

\begin{lemma}\label{lemme:boundrj}
Under Assumptions~\ref{ass:id} and \ref{ass:kernel}, we have the following upper-bound, for any positive integer $j$, 
\begin{equation*}
|r_j(\hat{T}_{h,\by}) - r_j(T_h)| \leq \| \hat{T}_{h,\by} - T_h\|_{HS},
\end{equation*}
where $r_j$ is defined by \eqref{eq:defrj}.
\end{lemma}

\begin{lemma} \label{lemma:eqnorms}
Let $\phi: \mathbb{R}^d\times\mathbb{R}^d \to \mathbb{R}^+$ be a square-integrable function with $$\|\phi\|_2 = \left[ \int \phi^2(\bz_1,\bz_2) \, \mathrm{d}\bz_1 \, \mathrm{d}\bz_2\right]^{1/2} < \infty$$
 and $\mathcal{T}:\spa\to\spa$ be the integral operator defined by
 $$
 [\mathcal{T}(\omega)](\bz_{2}) = \int_{\mathbb{R}^d} \omega(\bz_{1}) \phi(\bz_1,\bz_{2}) \, \mathrm{d}\bz_{1},
 $$
 then we have the following equivalence between the norms
 $
 \|\mathcal{T}\|_{HS} = \|\phi\|_2.
 $
\end{lemma}

\begin{lemma} \label{lemma:boundfunctionG}
 With the notations of Section~\ref{sec:estim},  let
$g:\mathbb{R}^{d\times (n+1) } \to \mathbb{R}^+$ be given, for any $\bY=(\bY_1, \ldots, \bY_{n+1})$, by
$$
g(\bY)= \|\hat{p}_{h,\bY} - \mathbb{E}[K_h^{d}(\cdot - \bY_1)K_h^{d}(\cdot - \bY_2 )\|_2.
$$
Then, under Assumption~\ref{ass:kernel}, for any  $\by\in\mathbb{R}^{d(n+1)}$ and $\tilde\by\in\mathbb{R}^{d(n+1)}$,
$$
|g(\by) - g(\tilde\by)| \leq \frac{2\sqrt{2}}{n} \|K_h\|_2^{2d}\sum_{t=1}^{n+1} \mathds{1}_{\{\by_t\neq\tilde\by_t\}}.
$$
\end{lemma}

\begin{lemma} \label{lemma:boundexpectationG} Assume that Assumptions~\ref{ass:id} and \ref{ass:kernel} are fulfilled. Then, for any $h>0$, we have

\begin{equation}
\mathbb{E}[g(\bY)]\leq \frac{\|K_h\|^{2d}_{2}}{n^{1/2}} (9+8 t_{\text{mix}})^{1/2}
\end{equation}
with $t_{\text{mix}}$ defined in Lemma \ref{lemma: Mixing_Hidden} and $g$ defined in Lemma~\ref{lemma:boundfunctionG}.
\end{lemma}

\begin{lemma}[McDiarmid's inequality for Markov Chains \citep{paulin2015concentration}]\label{lemma:McDiarmid}

\begin{itemize}
\item Let $\bV=(V_1,\ldots,V_n)^\top$ be a (not necessarily time homogeneous) Markov chain, taking values in a Polish state space $\Lambda=\Lambda_1 \times \ldots \times \Lambda_n$. Suppose that $g:\Lambda \to \mathbb{R}$ is such that there exists some $\bc=(c_1,\ldots,c_n)$, which satisfies that, for any $\bv\in\Lambda$ and $\tilde \bv\in\Lambda$,
$$
|g(\bv) - g(\tilde{\bv})|\leq \sum_{t=1}^n c_t \mathds{1}_{\{v_t\neq \tilde{v}_t\}}.
$$
Then for any $t\geq 0$, we have
$$
\mathbb{P}\left(|g(\bV) - \mathbb{E}g(\bV)|\geq t \right)\leq 2 \exp \left( \frac{-2t^2}{9\|\bc\|^2 t_{\text{mix}}} \right),
$$
where
$t_{\text{mix}}$ is the mixing time of the Markov chain defined in Lemma \ref{lemma: Mixing_Hidden}.
\item Let $\bW=(W_1,\ldots,W_n)^\top$ be an hidden Markov chain with underlying chain $\bV=(V_1,\ldots,V_n)^\top$ having mixing time $t_{\text{mix}}$. Suppose that $g:\Lambda \to \mathbb{R}$ satisfies that for any $\bw\in\Lambda$ and $\tilde \bw\in\Lambda$
\begin{equation} \label{eq:boundeddiff}
|g(\bw) - g(\tilde{\bw})|\leq \sum_{t=1}^n c_t \mathds{1}_{\{w_t\neq \tilde{w}_t\}},
\end{equation}
for some $\bc=(c_1,\ldots,c_n)$, then for any $t\geq 0$, we have
$$
\mathbb{P}\left(|g(\bW) - \mathbb{E}g(\bW)|\geq t \right)\leq 2 \exp \left( \frac{-2t^2}{9\|\bc\|^2 t_{\text{mix}}} \right),
$$
\end{itemize}
\end{lemma}

\section{Proofs of the main results}

\begin{proof}[Proof of Lemma~\ref{lemma:ident}]
Since, by  Assumption~\ref{ass:id}.\ref{ass:emission}, the densities $\{f_1,\ldots,f_L\}$ are linearly independent, there exist $\kappa_1,\ldots,\kappa_L$ with $\kappa_\ell\in\mathbb{R}^d$ such that the $L\times L$ matrix $\bM_f$ defined by $\bM_f[\ell,j]=f_\ell(\kappa_j)$, has full rank.
Let the $L\times L$ matrix $\bM_g$ be defined by $\bM_g[\ell,j]=g_\ell(\kappa_j)$. From \eqref{eq:distcouple}, we have $\bM_g=\bA M_f$. Noting that $\text{det}(\bM_g)=\text{det}(\bA)\text{det}(\bM_f)$ and that $\bA$ is invertible by Assumption~\ref{ass:id}.\ref{ass:transition}, we have that $\bM_g$ has full rank and thus that the densities $\{g_1,\ldots,g_L\}$ are linearly independent. Thus, using Proposition 3 of \cite{KasaharaJRSSB2014}, we obtain that $L$ is identifiable from the distribution of a pair of consecutive observations.
\end{proof}

 \begin{proof}[Proof of Lemma~\ref{lemma:controlSVD}]
 By Lemma~\ref{lemma:Hoffman}, we have under Assumptions~\ref{ass:id} and \ref{ass:kernel}
\begin{equation} \label{eq:sigmaL}
\sum_{j=1}^L\left( \sigma_j(T_h) - \sigma_j(T)\right)^2 \leq \| T_h - T\|_{HS}^2.
\end{equation}
We now show that the right-hand side of \eqref{eq:sigmaL} is of order $h^4$. First, note that $T_h-T$ is an integral operator given, for any function $\omega \in \spa$, by
$$[T_h-T](\omega)(\bz_1) = \int_{\mathbb{R}^d} \omega(\bz_1) \left[ p_h(\bz_1, \bz_2) - p(\bz_1,\bz_2)\right] \mathrm{d} \bz_2, $$
with $p_h$ defined, for some bandwidth $h >0$ and some kernel $K_h^d(\bz)= \prod_{j=1}^d K_{hj}(z_j)$, by $p_h(\bz_1,\bz_2)=\int_{\mathbb{R}^{d}\times\mathbb{R}^{d}}p(\by_1,\by_2)K^d_h(\bz_1-\by_1)K^d_h(\bz_2-\by_2)\, \mathrm{d}\by_1 \, \mathrm{d}\by_2$.
From Lemma~\ref{lemma:eqnorms}, we thus have
$$
\|T_h - T \|^2_{HS} = \|p_h-p\|_2^2.$$
Besides, variable change theorem implies that
$$
p_h(\bz_1,\bz_2)= h^{2d} \int_{\mathbb{R}^{d}\times\mathbb{R}^{d}} p(\bz_1- h \bu_1 ,\bz_2- h \bu_2 ) K_h^d(h\bu_1) K_h^d(h\bu_2) \dd\bu_1 \dd\bu_2.
$$
Therefore, a Taylor expansion of order 2 of $p(\bz_1-\bu_1 h,\bz_2-\bu_2 h)$ around $(\bz_1^\top,\bz_2^\top)^\top$ and Assumptions~\ref{ass:regularity} implies that
\begin{equation}\label{eq:bias}
\|T_h - T \|^2_{HS} = O(h^4).
\end{equation}
 \end{proof}

\begin{proof}[Proof of Theorem~\ref{thm:finitesample}]
By  Lemma~\ref{lemme:boundrj}, applied with $j=L+1$, 
we have, for any positive $\tau$, the following inclusion of events
$$ \{\|\hat{T}_{h,\by} - T_h\|_{HS} \leq \tau \}\subseteq\{ \vert r_{L+1}(\hat{T}_{h,\by}) - r_{L+1}(T_h) \vert \leq \tau \} , $$
which leads to this keystone inclusion
\begin{equation}\label{eq:inclusionevents} 
\{\|\hat{T}_{h,\by} - T_h\|_{HS} \leq \tau \}\subseteq \{\hat{L}(\tau,h)\leq L \}.
\end{equation} 
Indeed, using the fact that $T_h$ is of rank $L$, we have $r_{L+1}(T_h)=0$, which implies that $$\{\hat{L}(\tau,h)\leq L \}= \{r_{L+1}(\hat{T}_{h,\by}) \leq \tau \} = \{ \vert r_{L+1}(\hat{T}_{h,\by}) - r_{L+1}(T_h) \vert \leq \tau \}.$$ 
Thus,  controlling the probability that $\hat{L}(\tau,h)$ overestimates $L$ can be achieved via a concentration inequality on $\|\hat{T}_{h,\by} - T_h\|_{HS}$. Noting that $\mathbb{E}\hat{T}_{h,\by}=T_h$,  Lemma~\ref{lemma:eqnorms} implies
$$
\|\hat{T}_{h,\by}- T_h\|_{HS} = g(\bY),
$$
with $g:\mathbb{R}^{d\times (n+1) } \to \mathbb{R}^+$ and 
$$
g(\bY)= \|\hat{p}_{h,\bY} - \mathbb{E}[K_h^{d}(\cdot - \bY_1)K_h^{d}(\cdot - \bY_2 )]\|_2.
$$
Therefore, a concentration inequality on $\| \hat{T}_{h,\by} - T_h\|_{HS}$ can be obtained from a concentration inequality on $g(\bY)$.  From Lemma~\ref{lemma:boundfunctionG}, we have  for any  $\by\in\mathbb{R}^{d(n+1)}$ and $\tilde\by\in\mathbb{R}^{d(n+1)}$, 
\begin{equation}
|g(\by) - g(\tilde\by)| \leq \frac{2\sqrt{2}}{n} \|K_h\|_2^{2d}\sum_{t=1}^{n+1} \mathds{1}_{\{\by_t\neq\tilde\by_t\}}.
\end{equation}
Thus, a concentration inequality on $g(\bY)$ can be achieved by a McDiarmid's inequality for Markov Chains (see Lemma~\ref{lemma:McDiarmid}) stated by \citet{paulin2015concentration} since condition \eqref{eq:boundeddiff} is satisfied.
Therefore, noting that vector $c$ defined in  Lemma~\ref{lemma:McDiarmid} is here a vector of length $n+1$ where each element is equal to $2\sqrt{2} \|K_h\|_2^{2d}/n$,  we have for any $t>0$
\begin{equation}\label{eq:concentre}
\mathbb{P}( \| \hat{T}_{h,\by} - T_h\|_{HS}\geq t + \mathbb{E}\left[ g(\bY) \right]) \leq \exp \left( -\frac{n^2}{n+1}\frac{t^2}{36 \|K_h\|_2^{4d} t_{\text{mix}}} \right).
\end{equation}
Thus, we can rewrite that, for any $0<\alpha<1$
\begin{equation}\label{eq:concentration}
\mathbb{P}\left( \| \hat{T}_{h,\by} - T_h\|_{\spa}\geq \left(C_{\alpha,1} \|K_h\|_2^{4d}(n+1)/n^2\right)^{1/2} + \mathbb{E}\left[ g(\bY) \right]\right)\leq \alpha ,
\end{equation}
where $C_{\alpha,1}=36 \ln(1/\alpha) t_{\text{mix}} $.
  Moreover, from Lemma~\ref{lemma:boundexpectationG}, we have
\begin{equation*}
\mathbb{E}[g(\bY)]\leq \bigg(\frac{\|K_h\|^{4d}_{2}}{n} C_2 \bigg)^{1/2},
\end{equation*}
where $C_2= 9 + 8 \, t_{\text{mix}}$.
Therefore, replacing $\mathbb{E}\left[ g(\bY) \right]$ by its upper-bound  in \eqref{eq:concentration}, we obtain that for any $0<\alpha<1$
\begin{equation}\label{eq:concentration2}
\mathbb{P}\left( \| \hat{T}_{h,\by} - T_h\|_{HS}\geq \frac{\|K_h\|_2^{2d}}{n^{1/2}} \left( \left(C_{\alpha,1} \frac{n+1}{n}\right)^{1/2} + C_2^{1/2} \right) \right)\leq \alpha .
\end{equation}
Combining \eqref{eq:inclusionevents} and \eqref{eq:concentration2} leads to \eqref{eq:overestimation_nknown}.

In addition, to obtain \eqref{eq:underestimation_nknown}, it is important to first notice the following equality of events
$$
\{\hat{L}(\tau,h) = L\} = \left\{ \{r_L( \hat{T}_{h,\by} ) > \tau \} \cap\{ r_{L+1}(\hat{T}_{h,\by}) < \tau \} \right\} .
$$
We then recall that $r_{L+1}(T_h)=0$ and that $r_{L}(T_h)=\sigma_L(T_h)$. Thus, by Lemma~\ref{lemme:boundrj} applied with $j=L+1$ and $j=L$ respectively, we obtain that $$r_{L+1}(\hat{T}_{h,\by})\leq \|\hat{T}_{h,\by} - T_h\|_{HS},$$
and $$ r_L(\hat{T}_{h,\by})\geq \sigma_L(T_h) - \|\hat{T}_{h,\by} - T_h\|_{HS}. $$ Therefore, on the event $\{\sigma_L(T_h)>2\tau\}\cap\{\|\hat{T}_{h,\by} - T_h\|_{HS}<\tau\}$, we have $r_{L+1}(\hat{T}_{h,\by})<\tau$ and   $r_L(\hat{T}_{h,\by})\geq \tau$, which leads to the following inclusion of events
$$
\left\{ \{\sigma_L(T_h)>2\tau\} \cap \{\|\hat{T}_{h,\by} - T_h\|_{HS}<\tau\}\right\}\subseteq \{\hat{L}(\tau,h) = L\}.
$$
For any $0<\alpha<1$, the event $\{\sigma_L(T_h)>2\tau_{\alpha}\} $ is not random (see \eqref{eq:tau} for the definition of $\tau_{\alpha}$). Thus, if there exists an $\tilde h$ such that $\sigma_L(T_{\tilde{h}})>2\tau_{\alpha}$, we can then conclude that 
$$
\mathbb{P}(L(\tau_{\alpha},\tilde h)=L)\geq 1 - \alpha.
$$
To complete the proof, we have to show that such an $\tilde{h}$ exists. Note that using Lemma~\ref{lemma:controlSVD} and the assumption that $\sigma_L(T)>2\tau_{\alpha} + \varepsilon$, we obtain that
$$
\sigma_L(T_h)>2\tau_{\alpha} + \varepsilon + O(h^2),
$$
which ensures the existence of $\tilde h$.
\end{proof}

\begin{proof}[Proof of Proposition~\ref{prop:Kmeans}]
Let $\mu_{a,b}=\frac{1}{b-a+1}\sum_{\ell=a}^b \wwr_\ell$, we have 
$$\mu_{1,s+1}=\frac{s}{s+1} \mu_{1,s} + \frac{1}{s+1} \wwr_{s+1}$$
and
$$\mu_{s+1,n}=\frac{n-s-2}{n-s-1} \mu_{s+2,n} + \frac{1}{n-s-1} \wwr_{s+1}.$$
Thus, we have for any $s\in\{1,\ldots,m-1\}$
\begin{align*}
    \sum_{\ell=1}^{s+1}(\wwr_\ell - \mu_{1,s+1})^2 &=\sum_{\ell=1}^s (\wwr_\ell - \mu_{1,s+1})^2 + s(\mu_{1,s} - \mu_{1,s+1})^2 + (\wwr_{s+1} - \mu_{1,s+1})^2 \\
    &= \sum_{\ell=1}^s (\wwr_\ell - \mu_{1,s})^2 + \frac{s}{s+1} (\mu_{1,s} - \wwr_{s+1})^2,
\end{align*}
and
\begin{align*}
    \sum_{\ell=s+1}^{m}(\wwr_\ell - \mu_{s+1,m})^2 &=\sum_{\ell=s+2}^m (\wwr_\ell - \mu_{s+2,m})^2 + (n-s-2)(\mu_{s+2,m} - \mu_{s+1,m})^2 + (\wwr_{s+1} - \mu_{s+2,m})^2 \\
    &= \sum_{\ell=s+2}^m (\wwr_\ell - \mu_{s+2,m})^2 + \frac{m-s-2}{m-s-1} (\wwr_{s+1} - \mu_{s+2,m})^2.
\end{align*}
The K-means algorithm aims to minimize the within-group variance. Since the $\wwr_\ell$ are ordered, by considering two groups, the K-means algorithm aims to find $s^\star$ such that
$$
s^\star = \argmin_{s\in\{2,\ldots,m-1\}} \Delta(s),
$$
where
$$
\Delta(s) = \sum_{\ell=1}^s(\wwr_\ell - \mu_{1,s})^2 + \sum_{\ell=s+1}^m (\wwr_\ell - \mu_{s+1,m})^2.
$$
Therefore, the proof is complete if we show that $s^\star=L$ and so that $\Delta(s)>\Delta(L)$ for any $s\neq L$.
For any $s<L-1$, we have
$$
\Delta(s) - \Delta(s+1) = \frac{m-s-2}{m-s-1} (\wwr_{s+1} - \mu_{s+2,m})^2 - \frac{s}{s+1} (\mu_{1,s} - \wwr_{s+1})^2.
$$
Noting that $\mu_{s+2,m}< \wwr_{s+2} (L-s-2)/(n-s-2) + \varepsilon (n-L)/(n-s-2) \leq \wwr_{s+1} (L-s-2)/(n-s-2) + \varepsilon (n-L)/(n-s-2)$, then
$$
\frac{m-s-2}{m-s-1} (\wwr_{s+1} - \mu_{s+2,m})^2 > \frac{(n-L)^2}{(n-s-1)(n-s-2)} (\wwr_{s+1} - \varepsilon)^2.
$$
Since $\mu_{1s}< [ (s/\xi (s+1))^{1/2} + 1]\wwr_{s+1}$ for $s<L$, and since $\mu_{1s}>\wwr_{s+1}$, we have
$$
\frac{s}{s+1} (\mu_{1,s} - \wwr_{s+1})^2>\frac{1}{\xi}\wwr_{s+1}^2.
$$
Therefore, 
$$
\Delta(s) - \Delta(s+1) > \frac{\xi - 1}{\xi}\wwr_{s+1}^2 + o(m) + O(\varepsilon).
$$
and thus since $\xi>1$, if $m$ is large enough and if $\varepsilon$ is small enough, then for $s\in\{1,\ldots,L-1\}$, we have
$\Delta(s) - \Delta(s+1) >0$ implying that for any $s\in\{1,\ldots,L-1\}$, we have $\Delta(s)>\Delta(L)$.

If $s>L+1$, then 
$$
\Delta(s)=\sum_{\ell=1}^L (\wwr_\ell - \mu_{1L})^2 + L(\mu_{1L} - \mu_{1s})^2 + (s-L)(\varepsilon - \mu_{1s})^2 + (m-s-1)O(\varepsilon^2).
$$
Thus, for any $s\in\{L+1,\ldots,m\}$
$$
\Delta(s) - \Delta(L) = (L-s-1) O(\varepsilon^2) + L(\mu_{1L} - \mu_{1s})^2 + (s-L)\mu_{1s}^2 + (s-L)O(\varepsilon).
$$
Therefore, if $\varepsilon$ is small enough  then $\Delta(s)>\Delta(L)$ for $s>L$.
\end{proof}

\section{Proofs of the preliminary lemmas}\label{App:proofl}

\begin{proof}[Proof of Lemma~\ref{lemma: Mixing_Hidden}]

Step 1: The Markov and stationary properties of the process $\vec{\bX}=(\vec{\bX}_1, \ldots, \vec{\bX}_n)$ is straightforward. It remains just to prove that its mixing time is controlled by the mixing time $t_{mix}$ of the process $(X_t)_{t\in \mathbb{Z}}$. We denote by ${A}^{\vec{\bX}}:\{1,\ldots, L\}^2\times \{1,\ldots, L\}^2$ the kernel transition and by $\bpi^{\vec{\bX}}$ its stationary distribution, which are given for all $(i,j,k,\ell) \in \{1,\ldots, L\}^4$ by
$$
{A}^{\vec{\bX}}((i,j),(k,\ell))=\mathbf{1}_{\{k=j\}}A[k,\ell] \, \; \text{ and } \, \;
\bpi^{\vec{\bX}}(i,j)=\mathbb{P}_{\bpi}(X_s=i, X_{s+1}=j)
.$$
Let $t\geq t_{mix}+1$ and $\delta_{(i_1,j_1)}$ the Dirac distribution on $\{1,\ldots, L\}^2$ that puts mass $1$ at the pair $(i_1,j_1)$ and $0$ everywhere else, the total variation distance between $(A^{\vec{\bX}})^t$ and the stationary distribution $\bpi^{\vec{\bX}}$ is then given by

\begin{align*}
\|\delta_{(i_1,j_1)}(A^{\vec{\bX}})^t-\bpi^{\vec{\bX}}\|_{TV}&=\frac{1}{2}\sum_{(i,j) \in \{1,\ldots, L\}^2} \vert \delta_{(i_1,j_1)}(A^{\vec{\bX}})^t[i,j]-\bpi^{\vec{\bX}}(i,j)\vert \\
&=\frac{1}{2}\sum_{(i,j) \in \{1,\ldots, L\}^2} \vert \mathbb{P}(\vec{\bX}_{t+1}=(i,j)\vert \vec{\bX}_{1}=(i_1,j_1))-A[i,j]\bpi(i)\vert\\
&=\frac{1}{2}\sum_{(i,j) \in \{1,\ldots, L\}^2} \vert \mathbb{P}((X_{t+1},X_{t+2})=(i,j)\vert (X_1,X_2)=(i_1,j_1))-A[i,j]\bpi(i)\vert\\
&=\frac{1}{2}\sum_{(i,j) \in \{1,\ldots, L\}^2} \vert \mathbb{P}(X_{t+2}=j \vert X_{t+1}=i)\mathbb{P}(X_{t+1}=i\vert X_{2}=j_1 )-A[i,j]\bpi(i)\vert\\
&=\frac{1}{2}\sum_{i \in \{1,\ldots, L\}} \left(\sum_{j \in \{1,\ldots, L\}} A[i,j] \right) \vert A^{t-1}[j_1,i]-\bpi(i)\vert\\
&=\frac{1}{2}\sum_{i \in \{1,\ldots, L\}}  \vert A^{t-1}[j_1,i]-\bpi(i)\vert \\
&=\|\delta_{j_1}A^{t-1}-\bpi\|_{TV}.
\end{align*}
By making use of Assumption~\ref{ass:id} and  the definition of the mixing time of $(X_t)_t$, we deduce that 
$$\|\delta_{(i_1,j_1)}(A^{\vec{\bX}})^t-\bpi^{\vec{\bX}}\|_{TV} \leq 1/4$$

Step 2: : The Markov property and the mixing rate of the hidden state $\vec{\bZ}$ is ensured from the stability properties of hidden chains, i.e. by taking the same strategy as in Step 1 for the Markov kernel $A^{\vec{\bZ}}$ and let $\delta_{(i_1,j_1,\tilde{\by}_1,\tilde{\by}_2)}$ the Dirac distribution on $\{1,\ldots, L\}^2\times \mathbb{R}^{2d}$ that puts mass $1$ at the pair $(i_1,j_1,\tilde{\by}_1,\tilde{\by}_2)$ and $0$ everywhere else, we have that

\begin{align*}
\|\delta_{(i_1,j_1,\tilde{\by}_1,\tilde{\by}_2)}(A^{\vec{\bZ}})^t-\bpi^{\vec{\bZ}}\|_{TV}
=\|\delta_{(j_1,\tilde{y}_2)}(A^{X,\bY})^{t-1}-\bpi^{X,\bY}\|_{TV} \leq 1/4
\end{align*}
where $A^{X,\bY}$ is the Markov kernel transition of the hidden chain $(X_t,\bY_t)$ and $\bpi^{X,\bY}$ its stationary distribution which is also uniformly ergodic by stability and Assumption~\ref{ass:id}.

\end{proof}

\begin{proof}[Proof of Lemma~\ref{lemma: Cov_Paulin_hidden}]
The proof of the first part is a direct consequence of Theorem 3.4 and Proposition 3.4 in \cite{paulin2015concentration} for the hidden Markov chains $\vec{\bZ} \sim (A^{\vec{\bZ}},\bpi^{\vec{\bZ}})$ defined in Lemma \ref{lemma: Mixing_Hidden} with mixing time $t_{\text{mix}}^{\vec{\bZ}}$.  The second part follows straightforwardly by taking for $\phi$ the following particular function
\begin{displaymath}
\phi :
\left\lbrace
  \begin{array}{rcl}
    \{1,\ldots, L \}^{2} \times \mathbb{R}^{2d} & \longrightarrow &\mathbb{R} \\
    (k, \ell, \vec{\by}_1) & \longmapsto & \psi(\vec{\by}_1) \\
  \end{array},
\right.
\end{displaymath}
with $\psi$ being a measurable function from $\mathbb{R}^{2d} \longrightarrow \mathbb{R}$. We get 
$$\sum_{t=1}^{n}\mathbb{E}[\psi(\vec{\bY}_1)\psi(\vec{\bY}_t)]\leq \frac{2}{\gamma_{ps}} \mathbb{V}[\psi(\vec{\bY}_1)],$$
where $\gamma_{ps}$ is the pseudo-spectral gap associated with the kernel transition $A^{\vec{\bZ}}$ (see  \cite{paulin2015concentration} for a definition). We conclude the proof by noting that, from Assumption~\ref{ass:id}, we have $\frac{2}{\gamma_{ps}}\leq 4 t_{\text{mix}}^{\vec{\bZ}}$.
\end{proof}

\begin{proof}[Proof of Lemma~\ref{lemme:boundrj}]
Using the reverse triangle inequality, we can easily show that, for any positive integer $j$
$$
|r_j(\hat{T}_{h,\by}) - r_j(T_h)| \leq \left[\sum_{i\geq j} \left( \sigma_i(\hat{T}_{h,\by}) - \sigma_i(T_h)\right)^2\right]^{1/2}.
$$
We conclude the proof by using Lemma~\ref{lemma:Hoffman},  which leads to the announced result
\begin{equation*} 
|r_j(\hat{T}_{h,\by}) - r_j(T_h)| \leq \| \hat{T}_{h,\by} - T_h\|_{HS}.
\end{equation*}
\end{proof}

 \begin{proof}[Proof of Lemma~\ref{lemma:eqnorms}]
Considering an orthonormal basis $\{e_k\}_{k=1}^{\infty}$ of $\spa$, we have
$$
\|\mathcal{T}\|_{HS}^2 = \sum_{k=1}^\infty \int e_k(\by_1) \phi(\by_1,\by_2) \dd\by_2 \dd\by_1.
$$
Therefore, denoting $\phi_{\by_2}=\phi(\cdot,\by_2)$, we have
$$
\|\mathcal{T}\|_{HS}^2 = \int\sum_{k=1}^\infty  <\phi_{\by_2}, e_k> \dd\by_2.
$$
Applying Parseval's inequality, we obtain
\begin{equation*}
\|\mathcal{T}\|_{HS}^2 = \int \|\phi_{\by_2}\|_2^2 \dd\by_2
 =\|\phi\|_2^2 ,
\end{equation*}
which concludes the proof.
\end{proof}

\begin{proof}[Proof of Lemma~\ref{lemma:boundfunctionG}]
For any $\by\in\mathbb{R}^{d(n+1)}$ and $\tilde\by\in\mathbb{R}^{d(n+1)}$, define $\Delta=\{t:\; \by_t\neq\tilde\by_{t}\}$. Let the vectors $\by^{(s)}$ defined by $\by^{(0)}=\by$ and for $s=1,\ldots,\text{card}(\Delta)$ by
$$
\by_t^{(s)} = \left\{ \begin{array}{rl}
\by_t^{(s-1)} & \text{if } t\neq t^{(s)}\\
\tilde \by_t & \text{if } t= t^{(s)}\\
\end{array}
\right.,$$
where $t^{(s)}$ denotes the element $s$ of $\Delta$. Applying the reverse triangle inequality, we have
$$
|g(\by^{(s-1)}) - g(\by^{(s)})| \leq \| \hat{p}_{h,\by^{(s-1)}} - \hat{p}_{h,\by^{(s)}}||_2.
$$
For any $\bz_1$ and $\bz_2$ in $\mathbb{R}^d$, using \eqref{eq:phat}, we have
$$\left[\hat{p}_{h,\by^{(s-1)}}- \hat{p}_{h,\by^{(s)}}\right](\bz_1,\bz_2) = \frac{1}{n} \sum_{t=1}^{n} \left( K_h^d(\bz_1-\by_t^{(s-1)}) K_h^d(\bz_2-\by_{t+1}^{(s-1)}) - K_h^d(\bz_1-\by_t^{(s)})K_h^d(\bz_2-\by_{t+1}^{(s)})\right). $$
Using the fact that, by construction, for any $t \neq t^{(s)}$, $\by_t^{(s)} = \by_t^{(s-1)}$, we deduce that all the terms in the previous sum vanish, except the ones for $t=t^{(s)}$ or $t=t^{(s)}-1$.
Thus, if $1<t^{(s)}<n$, we obtain that
\begin{multline*}
\| \hat{p}_{h,\by^{(s-1)}} - \hat{p}_{h,\by^{(s)}}||_2^2 =\frac{1}{n^2} 
\int
\left[ K_h^d( \bz_1-\by^{(s-1)}_{t^{(s)}-1}) \left( K_h^d( \bz_2-\by^{(s-1)}_{t^{(s)}} ) - K_h^d(\bz_2-\by^{(s)}_{t^{(s)}} ) \right)\right.\\
+ \left.
 \left( K_h^d(\bz_1-\by^{(s-1)}_{t^{(s)}}) - K_h^d(\bz_1-\by^{(s)}_{t^{(s)}} )\right) K_h^d(\bz_2 - \by^{(s-1)}_{t^{(s)}+1}) \right]^2
\dd\bz_1 \dd\bz_2.
\end{multline*}
Note that if $t^{(s)}=1$ or $t^{(s)}=n$, then the same reasoning can be applied but only one term appears in the integrand on the right-hand side of the previous equation. By Assumption~\ref{ass:kernel}, we have  that $\int
\left(K_h^d(\bz_1-\ba)K_h^d(\bz_2-\bb) \right)^2 \dd\bz_1 \dd\bz_2 = \|K_h\|_2^{4d}$ for any $\ba\in\mathbb{R}^d$ and $\bb\in\mathbb{R}^d$. Thus, using that $(a+b)^2\leq 4 (a^2 \wedge b^2)$ and  $(a-b)^2\leq 2 (a^2 \wedge b^2)$  for any $a>0$ and $b>0$, we deduce that
$$
|g(\by^{(s-1)}) - g(\by^{(s)})| \leq  \frac{2\sqrt{2}}{n}\|K_h\|_2^{2d}.
$$
The proof is completed by noticing that $\by^{(\text{card}(\Delta))}=\tilde\by$ and that 
$$
|g(\by) - g(\tilde\by)| \leq \sum_{s=1}^{\text{card}(\Delta)}  |g(\by^{(s-1)}) - g(\by^{(s)})|.
$$
\end{proof}

\begin{proof}[Proof of Lemma~\ref{lemma:boundexpectationG}]

Let us denote be $U_{t,n}(\bz)$ the random variable defined as $$U_{t,n}(\bz)=K_h^d(\bz_1-\bY_t)K_h^d(\bz_2-\bY_{t+1})$$ where $\bz=(\bz_1,\bz_2)$.  We recall that the function $g:\mathbb{R}^{d\times (n+1) } \to \mathbb{R}^+$ is such that $$g(\bY)= \|\hat{p}_{h,\bY} - \mathbb{E}[K_h^{d}(\cdot - \bY_1)K_h^{d}(\cdot - \bY_2 )\|_2$$ with
\begin{equation*}
\hat{p}_{h,\by}(\bz) = \frac{1}{n} \sum_{t=1}^{n} K_h^d(\bz_1-\by_t)K_h^d(\bz_2-\by_{t+1}). 
\end{equation*}
We also denote by $p_h(\bz)=\mathbb{E}[\hat{p}_{h,\by}(\bz)]$. 
Hence, using the concavity of the square-root function and the Jensen's inequality, we get
\begin{align*}
\mathbb{E}[g(\bY)]&=\mathbb{E}\left[\|\frac{1}{n}\sum_{t=1}^{n} U_{t,n}(\bz)-\mathbb{E}\left[U_{t,n}(\bz)\right] \|_2\right]\\
&\leq \mathbb{E}^{1/2}\left[\|\frac{1}{n}\sum_{t=1}^{n} (U_{t,n}(\bz)-\mathbb{E}[U_{t,n}(\bz)])\|_2^2\right]\\
&\leq \mathbb{E}^{1/2}\left[\|\hat{p}_{h,\by}(\bz)-p_{h}(\bz)\|_2^2\right]
\end{align*}
Using the standard bias variance decomposition, we obtain that
\begin{align*}
\mathbb{E}\left[\|\hat{p}_{h,\by}(\bz)-p_{h}(\bz)\|_2^2\right]&=\mathbb{E}\bigg[\int (\hat{p}_{h,\by}(\bz)-p_{h}(\bz))^2 \dd \bz\bigg]\\
&= \int \mathbb{E}\left[(\hat{p}_{h,\by}(\bz)-p_{h}(\bz))^2\right] \dd \bz \\
&=\int \mathbb{V}[\hat{p}_{h,\by}] \dd\bz \\
& =\int \mathbb{V}\left[\frac{1}{n}\sum_{t=1}^{n}U_{t,n}(\bz)\right] \dd\bz\\
&=\int\bigg(\frac{1}{n^2}\sum_{t=1}^{n}\mathbb{V}\left[U_{t,n}(\bz)\right]+\frac{2}{n^2}\sum_{1\leq t\leq t'\leq n}Cov(U_{t,n}(\bz), U_{t',n}(\bz))\bigg)\dd\bz\\
&=\int\bigg(\frac{1}{n}\mathbb{V}[U_{1,n}(\bz)]+\frac{2}{n^2}\sum_{t=2}^{n}(n-t) Cov(U_{1,n}(\bz), U_{t,n}(\bz))\bigg)\dd \bz\\ 
&\leq \int\bigg(\frac{1}{n}\mathbb{V}[U_{1,n}(\bz)]+\frac{2}{n}\sum_{t=2}^{n}Cov(U_{1,n}(\bz), U_{t,n}(\bz))\bigg)\dd \bz
\end{align*}
We start by computing the first term on the right-hand side 
\begin{align}
\frac{1}{n}\int \mathbb{V}[U_{1,n}(\bz)]d\bz&=\frac{1}{n}\int \mathbb{V}\bigg[K_h^d(\bz_1-\bY_1)K_h^d(\bz_2-\bY_2)\bigg]\dd \bz\notag\\
&\leq \frac{1}{n}\int \bigg(\int \int K_h^{2d}(\bz_1-\by_1)K_h^{2d}(\bz_2-\by_2)p(\by_1,\by_2) \dd\by_{1} \dd\by_{2}\bigg) \dd \bz\notag\\
&\leq \frac{1}{n}\int \int \int K_h^{2d}(\mathbf{u})K_h^{2d}(\mathbf{v})p(\bz_1- \mathbf{u},\bz_2-\mathbf{v}) \dd \mathbf{u} \dd \mathbf{v} \dd \bz\notag\\
&\leq \frac{1}{n}\|K_h\|^{4d}_{2}\label{eq:var}
\end{align}
Concerning the covariance terms, taking for $\psi$:
\begin{displaymath}
\psi :
\left\lbrace
  \begin{array}{rcl}
  \mathbb{R}^{2d} & \longrightarrow &\mathbb{R} \\
    \vec{\by}_1 & \longmapsto & \psi(\vec{\by}_1)= \displaystyle K_h^d(\bz_1-\by_1)K_h^d(\bz_2-\by_2)\\
  \end{array}
\right.
\end{displaymath}
we have from Lemma \ref{lemma: Cov_Paulin_hidden}
\begin{align*}
&\sum_{t=2}^{n} Cov(U_{1,n}(\bz),U_{t,n}(\bz))\\
&=\sum_{t=2}^{n} Cov\bigg(K_h^d(\bz_1-\bY_1)K_h^d(\bz_2-\bY_2),K_h^d(\bz_1-\bY_t) K_h^d(\bz_2-\bY_{t+1})\bigg)\\
&\leq \sum_{t=1}^{n}\bigg(\mathbb{E}\bigg[\vert K_h^d(\bz_1-\bY_1)K_h^d(\bz_2-\bY_2)K_h^d(\bz_1-\bY_t)K_h^d(\bz_2-\bY_{t+1})\vert\bigg] \\
& \; \hspace{2cm} +\bigg(\mathbb{E}\left[\vert K_h^d(\bz_1-\bY_1)K_h^d(\bz_2-\bY_2)\vert\right]\bigg)^2\bigg) \\
&\leq 4 t_{mix}^{\vec{\bZ}} \; \mathbb{V}\bigg[ K_h^{d}(\bz_1-\bY_1)K_h^{d}(\bz_2-\bY_2)\bigg] \\
&\leq 4 t_{mix}^{\vec{\bZ}} \; \mathbb{E}\bigg[ K_h^{2d}(\bz_1-\bY_1)K_h^{2d}(\bz_2-\bY_2)\bigg]\\
&\leq 4 t_{mix}^{\vec{\bZ}} \bigg(\int \int K_h^{2d}(\bz_1-\by_1)K_h^{2d}(\bz_2-\by_2)p(\by_1,\by_2) \dd \by_{1} \dd \by_{2}\bigg),
\end{align*}
with $t_{\text{mix}}^{\vec{\bZ}}$ defined in Lemma \ref{lemma: Mixing_Hidden}. Integrating out the previous inequality over variable $\bz$, we obtain that
\begin{equation}\label{eq:cov:Paulin}
\int \sum_{t=2}^{n} Cov(U_{1,n}(\bz),U_{t,n}(\bz)) \dd\bz\leq 4 t_{\text{mix}}^{\vec{\bZ}}\|K_h\|^{4d}_{2}.
\end{equation}
 Combining \eqref{eq:var} and \eqref{eq:cov:Paulin}, we obtain \begin{align*}
\mathbb{E}[\|\hat{p}_{h,\by}(\bz)-p_{h}(\bz)\|_2^2]\leq \frac{\|K_h\|^{4d}_{2}}{n}(1+8t_{\text{mix}}^{\vec{\bZ}}).
\end{align*}
And, we can conclude with the announced bound, by making use of Lemma \ref{lemma: Mixing_Hidden}
\begin{equation*}
\mathbb{E}[g(\bY)]\leq \frac{\|K_h\|^{2d}_{2}}{n^{1/2} }(1+8t_{mix}^{\vec{\bZ}})^{1/2} \leq \frac{\|K_h\|^{2d}_{2}}{n^{1/2} }(9+ 8t_{\text{mix}})^{1/2}
\end{equation*}
\end{proof}

\section{Additionnal numerical experiments}

\begin{table}[htp!]
\centering
\begin{scriptsize}
\begin{tabular}{cc|cc  c c|cc c c|cc c c|cc  c c}
  \hline
Method & $n$ & \multicolumn{4}{c|}{Gaussian} & \multicolumn{4}{c|}{Student} & \multicolumn{4}{c|}{Laplace} & \multicolumn{4}{c}{Von-Mises}\\
& & L-1 & L-2 & \normalfont{L-3} & L$>$3 & L-1 & L-2 & \normalfont{L-3} & L$>$3 & L-1 & L-2 & \normalfont{L-3} & L$>$3 & L-1 & L-2 & \normalfont{L-3} & L$>$3 \\ 
  \hline
proposed & 250 & 0 & 19 & 81 & 0 & 0 & 9 & 91 & 0 & 0 & 9 & 91 & 0 & 0 & 5 & 84 & 11 \\ 
    & 500 & 0 & 2 & 98 & 0 & 0 & 2 & 98 & 0 & 0 & 1 & 99 & 0 & 0 & 0 & 96 & 4 \\ 
    & 1000 & 0 & 0 & 100 & 0 & 0 & 0 & 100 & 0 & 0 & 0 & 100 & 0 & 0 & 0 & 99 & 1 \\ 
    & 2000 & 0 & 0 & 100 & 0 & 0 & 0 & 100 & 0 & 0 & 0 & 100 & 0 & 0 & 0 & 100 & 0 \\ 
    & 4000 & 0 & 0 & 100 & 0 & 0 & 0 & 100 & 0 & 0 & 0 & 100 & 0 & 0 & 0 & 100 & 0 \\ \hline
  spectral & 250 & 0 & 0 & 1 & 99 & 0 & 0 & 0 & 100 & 0 & 0 & 0 & 100 & 0 & 0 & 1 & 99 \\ 
 slope & 500 & 0 & 0 & 0 & 100 & 0 & 0 & 0 & 100 & 0 & 0 & 0 & 100 & 0 & 0 & 2 & 98 \\ 
   & 1000 & 0 & 0 & 2 & 98 & 0 & 0 & 0 & 100 & 0 & 0 & 0 & 100 & 0 & 0 & 2 & 98 \\ 
    & 2000 & 0 & 0 & 0 & 100 & 0 & 0 & 0 & 100 & 0 & 0 & 0 & 100 & 0 & 0 & 2 & 98 \\ 
   & 4000 & 0 & 0 & 4 & 96 & 0 & 0 & 0 & 100 & 0 & 0 & 0 & 100 & 0 & 0 & 2 & 98 \\ \hline
  spectral & 250 & 16 & 18 & 20 & 43 & 11 & 29 & 19 & 38 & 15 & 20 & 30 & 34 & 13 & 26 & 32 & 28 \\ 
  kmeans  & 500 & 15 & 24 & 44 & 16 & 14 & 28 & 39 & 19 & 18 & 29 & 39 & 14 & 11 & 31 & 40 & 16 \\ 
  & 1000 & 4 & 17 & 76 & 3 & 5 & 18 & 69 & 8 & 6 & 16 & 73 & 5 & 0 & 20 & 75 & 5 \\ 
   & 2000 & 0 & 0 & 100 & 0 & 1 & 3 & 94 & 2 & 0 & 4 & 96 & 0 & 0 & 1 & 98 & 1 \\ 
   & 4000 & 0 & 0 & 100 & 0 & 0 & 0 & 100 & 0 & 0 & 0 & 100 & 0 & 0 & 0 & 100 & 0 \\ 
   \hline
\end{tabular}
\end{scriptsize}
\caption{Percentage of number of states selected by the competing methods (proposed method "proposed", spectral method with slope heuristic used for tuning threshold "spectral slope" and spectral method with the proposed method based on Kmeans used for tuning the constant "spectral kmeans"), according to the family of the emission distributions and the sample size, obtained on 100 replicates generated with $d=1$, $\nu=0.1$, with an marginal overlap between the emission distributions of 2.5\%.  }\label{tab:simudelta0025}
\end{table}

\begin{table}[htp!]
\centering
\begin{scriptsize}
\begin{tabular}{cc|cc  c c|cc c c|cc  c c|cc  c c}
  \hline
Method & $n$ & \multicolumn{4}{c|}{Gaussian} & \multicolumn{4}{c|}{Student} & \multicolumn{4}{c|}{Laplace} & \multicolumn{4}{c}{Von-Mises}\\
& & L-1 & L-2 & \normalfont{L-3} & L$>$3 & L-1 & L-2 & \normalfont{L-3} & L$>$3 & L-1 & L-2 & \normalfont{L-3} & L$>$3 & L-1 & L-2 & \normalfont{L-3} & L$>$3 \\
\hline
proposed & 250 & 0 & 72 & 28 & 0 & 0 & 82 & 18 & 0 & 0 & 69 & 31 & 0 & 0 & 9 & 90 & 1 \\ 
    & 500 & 0 & 75 & 25 & 0 & 0 & 86 & 14 & 0 & 0 & 65 & 35 & 0 & 0 & 2 & 98 & 0 \\ 
    & 1000 & 0 & 51 & 49 & 0 & 0 & 31 & 69 & 0 & 0 & 17 & 83 & 0 & 0 & 0 & 100 & 0 \\ 
    & 2000 & 0 & 1 & 99 & 0 & 0 & 0 & 100 & 0 & 0 & 0 & 100 & 0 & 0 & 0 & 100 & 0 \\ 
    & 4000 & 0 & 0 & 100 & 0 & 0 & 0 & 100 & 0 & 0 & 0 & 100 & 0 & 0 & 0 & 100 & 0 \\ \hline
  spectral & 250 & 0 & 2 & 12 & 86 & 0 & 0 & 0 & 100 & 0 & 0 & 1 & 99 & 0 & 4 & 9 & 87 \\ 
   slope & 500 & 0 & 4 & 17 & 79 & 0 & 0 & 7 & 93 & 0 & 0 & 1 & 99 & 0 & 0 & 14 & 86 \\ 
    & 1000 & 0 & 0 & 32 & 68 & 0 & 0 & 4 & 96 & 0 & 0 & 2 & 98 & 0 & 0 & 22 & 78 \\ 
    & 2000 & 0 & 0 & 39 & 61 & 0 & 0 & 8 & 92 & 0 & 0 & 1 & 99 & 0 & 0 & 22 & 78 \\ 
    & 4000 & 0 & 0 & 47 & 53 & 0 & 0 & 7 & 93 & 0 & 0 & 1 & 99 & 0 & 0 & 29 & 71 \\  \hline
  spectral & 250 & 21 & 11 & 23 & 41 & 17 & 21 & 25 & 35 & 12 & 19 & 30 & 38 & 9 & 31 & 25 & 35 \\ 
  kmeans & 500 & 16 & 28 & 46 & 10 & 14 & 35 & 34 & 17 & 16 & 34 & 35 & 15 & 11 & 45 & 29 & 15 \\ 
    & 1000 & 5 & 24 & 69 & 2 & 7 & 24 & 62 & 7 & 6 & 32 & 60 & 2 & 2 & 19 & 75 & 4 \\ 
   & 2000 & 0 & 2 & 98 & 0 & 1 & 1 & 98 & 0 & 1 & 4 & 95 & 0 & 0 & 2 & 98 & 0 \\ 
    & 4000 & 0 & 0 & 100 & 0 & 0 & 0 & 100 & 0 & 0 & 0 & 100 & 0 & 0 & 0 & 100 & 0 \\ 
   \hline
\end{tabular}
\end{scriptsize}
\caption{Percentage of number of states selected by the competing methods (proposed method "proposed", spectral method with slope heuristic used for tuning threshold "spectral slope" and spectral method with the proposed method based on Kmeans used for tuning the constant "spectral kmeans"), according to the family of the emission distributions and the sample size, obtained on 100 replicates generated with $d=1$, $\nu=0.1$, with an marginal overlap between the emission distributions of 10\%. }\label{tab:simudelta01}
\end{table}

\end{document}